\documentclass[english,aip,cha,amsmath,amssymb,reprint]{revtex4-1}
\usepackage[T1]{fontenc}
\usepackage[latin9]{inputenc}
\usepackage{babel}
\usepackage{amsmath}
\usepackage{graphicx}
\usepackage{esint}
\usepackage[unicode=true]
 {hyperref}
\usepackage{color}
\usepackage{color, colortbl}
\definecolor{Gray}{gray}{0.9}
\makeatletter

\providecommand{\tabularnewline}{\\}

\usepackage{circuitikz}
\usepackage{cleveref}
\usepackage{mathtools}
\usepackage{graphicx}
\usepackage{hyperref}

\makeatother

\begin{document}

\title{Incorporating Inductances in Tissue-Scale Models of Cardiac Electrophysiology}

\author{Simone \surname{Rossi}}
\email[Email: ]{simone.rossi@unc.edu}
\affiliation{Department of Mathematics, University of North Carolina, Chapel
Hill, NC USA}

\author{Boyce E.~\surname{Griffith}}
\affiliation{Departments of Mathematics and Biomedical Engineering and McAllister Heart Institute, University of North Carolina, Chapel
Hill, NC USA}

\begin{abstract}
In standard models of cardiac electrophysiology, including the bidomain and monodomain models, local perturbations
can propagate at infinite speed.
We address this unrealistic property
by developing a hyperbolic bidomain model that is based on a generalization of Ohm's law with a
Cattaneo-type model for the fluxes.  Further,
we obtain a hyperbolic monodomain model in the case that the intracellular and extracellular conductivity tensors have the same anisotropy ratio.
In one spatial dimension, the hyperbolic
monodomain model is equivalent to a cable model that includes axial inductances, and
the relaxation times of the Cattaneo fluxes are strictly
related to these inductances.
 A purely linear analysis shows that
the inductances are negligible, but models of cardiac electrophysiology are highly nonlinear,
and linear predictions may not capture the fully nonlinear dynamics. In fact,
contrary to the linear analysis, we show that for simple nonlinear
ionic models, an increase in conduction velocity is obtained for small
and moderate values of the relaxation time. A similar behavior is
also demonstrated with biophysically detailed ionic models.
Using the Fenton-Karma model along with a low-order finite element spatial discretization, we
numerically analyze differences between the standard monodomain model and
the hyperbolic monodomain model.
In a simple benchmark test, we show that the propagation of the action potential is strongly influenced by the alignment of the fibers with respect to the mesh in both the parabolic and hyperbolic models when using relatively coarse spatial discretizations.
Accurate predictions of the conduction velocity require computational mesh spacings on the order of a single cardiac cell.
We also compare the two formulations in the case of spiral break up and atrial fibrillation in an anatomically detailed model of the left atrium, and we examine the effect of intracellular and extracellular inductances on the virtual electrode phenomenon.
\end{abstract}

\selectlanguage{english}%


\maketitle
\begin{quotation}

Since its introduction by \citet{Hodgkin}, the cable equation and its higher-dimensional generalizations \cite{Tung1978} have been common models of electrical impulse propagation in excitable media, including neurons and muscle.
The effects of inductances in these systems are considered to be relatively small, and so they are neglected in classical versions of these models.
By omitting such terms,
the standard equations of cardiac electrophysiology become parabolic, but,
as in all parabolic equations, local perturbations
can propagate at infinite speed. This unrealistic property
has been addressed in models of neurons by \citet{lieberstein1970source}, and
hyperbolic models including inductances have been proposed by \citet{Liebeastein,Engelbrecht1981}, and \citet{Engelbrecht}.
These models are also supported by the fact that neurons, skeletal
muscle cells, and cardiomyocytes show the typical resonance effects
due to inductances, as demonstrated by the studies of \citet{Clapham1982} and \citet{Koch1984}.
The common conclusion that inductances are negligible, which is based on the
linear analyses for neurons by Scott \cite{Scott, Scott1972}
and the numerical results of \citet{Kaplan1970}, may not be valid
in the complex arrangement of cardiac tissue, where the inhomogeneities
together with highly nonlinear reactions can lead to reentrant waves
and chaotic behavior. We thus derive a hyperbolic model
for cardiac electrophysiology, and we compare the solutions with
parabolic models in several cases, including simple excitation patterns as well
as for spiral waves, atrial fibrillation, and the virtual electrode phenomenon.

\end{quotation}

The heart is a complex organ with a highly heterogeneous structure.
Muscle cell are embedded in an extracellular compartment that includes many components, including
capillaries, connective tissue, and collagen.
The structural arrangement of the tissue is known to influence electrical impulse
propagation \cite{Spach1983,Vera1995,Bueno-Orovio2008}.
The standard models of electrophysiology neglect all these complexities,
and the resulting equations have the unrealistic feature that local
perturbations can propagate infinitely fast. These models are fundamentally
based on the cable equation. The works of \citet{Liebeastein} and \citet{lieberstein1970source}
were the first to suggest that the cable model for neurons should
contain inductances because of the three-dimensional nature of the
axon. Several authors, including \citet{Kaplan1970,Scott}, and \citet{Engelbrecht1981}, have
investigated this hypothesis for nerves. Their conclusions
were that inductances in neurons are of the order of few $\mu$H and
therefore are negligible, following a linear analysis
of a one-dimensional nerve discussed by \citet{Scott1972}. Several
years later, further studies by \citet{DeHaanRobertLandDeFelice1978,Clapham1982}, and \citet{Koch1984}
showed that the cell membranes of excitable tissue exhibit
self-inductance. In particular, \citet{Clapham1982} showed that the impedance of
the embryonic heart cell membrane resonates at a frequency
around 1 Hz, thereby enhancing homogeneity of the voltage. To our knowledge,
there has been no experimental study that addresses the role of inductances
for reentrant waves. In this chaotic scenario, inductances may have
an important role in the system dynamics. 

Excitable tissues have a characteristic speed of transmission that limits
the velocity at which signals can propagate. Any signal, including action potential
wavefronts, cannot propagate faster than this characteristic speed.
On the other hand, it is important to notice that propagation
of the action potential is related to the nonlinear dynamics of the
system, and not to the wave-like behavior that may be induced by inductances.
One way to transform the cable equation into a hyperbolic equation
is to simply add an ``inertial'' term proportional to the persistence
time of the diffusive process, as done by \citet{Zemskov2015} Unfortunately,
this type of model cannot be derived from reasonable physical assumptions. However, as we show here, it is possible to
perform a phenomenological derivation of a hyperbolic reaction-diffusion model
by using the Cattaneo model for the fluxes. This formulation was originally introduced
by \citet{cattaneo48} to eliminate the anomalies found using Fourier's
law in the heat equation, and this model has subsequently been used
in a wide range of applications, including forest fire models \cite{MeNdez},
chemical systems \cite{Gorecki1995,Lemarchand1998}, thermal combustion
\cite{Fort2004}, and the spread of viral infections \cite{Fort2002}.
Cattaneo-type models for the fluxes have been derived
in several ways, ranging from phenomenological and thermodynamical
derivations to isotropic and anisotropic random walks with reactions
\cite{Jou1996,MaEndez1998}. 
The use of Cattaneo-type
fluxes in a monodomain model of cardiac electrophysiology leads to
two additional terms in the equations proportional to the characteristic
relaxation time of the medium: the second derivative in time of the
potential, which is associated with ``inertia'', and the time derivative of the
ionic currents. Even if the relaxation time is small, the rapid variation
of the ionic currents can impart a contribution that may not be negligible.
This is particularly relevant near the front of the wave, where
fast currents give rise to the upstroke of the action potential. 


Verification and validation remain major challenges in computational electrophysiology.
\citet{Krishnamoorthi2014} propose that
the ``wave speed should not be sensitive to choices of numerical
solution protocol, such as mesh density, numerical integration scheme,
etc.'' Although this is a fundamentally desirable criterion, it is also nearly
impossible to achieve in practice. What does seem reasonable is to require that that the error
in the wave speed be ``small'', but how small the error
must be taken clearly depends on the case under consideration. For example,
for a single heart beat, an error of 5\% might represent a reasonable approximation.
On the other hand, models of cardiac fibrillation may require a
much smaller error. In fact, in this case, a 5\% error in the conduction
velocity can determine whether reentrant waves form, or when and how they
break up. Although this problem has been discussed in detail in prior work \cite{Pezzuto2016},
the common belief within the field seems to be that spatial discretizations on the order of 200~$\mu$m
are sufficient to capture the conduction velocities of the
propagating fronts. This estimate seems to hold for isotropic propagation at normal coupling strengths, but
in many important cases, the most relevant propagation of the electrical
signal is transverse to the alignment of the cardiac cells, where
the conductivities are typically 8 times smaller than in the longitudinal direction \cite{hand2009deriving}.
As shown by \citet{Quarteroni2017}, the mesh sizes needed to resolve
transverse propagation are actually closer to 25~$\mu$m. Such a small
mesh spacing requires the use of large-scale simulations and highly efficient
codes. Further, the need to use such high spatial resolutions challenges the fundamental idea of
using a continuum model for the description of cardiac electrophysiology,
and multiscale models have been proposed \cite{hand2010adaptive}.
This paper shows that the transverse conduction velocities are sensitive to the grid size and to the mesh orientation for both regular and hyperbolic versions of the model, and that the hyperbolic model has similar mesh size requirements as the standard model.

\begin{figure*}
\centering{}
\includegraphics{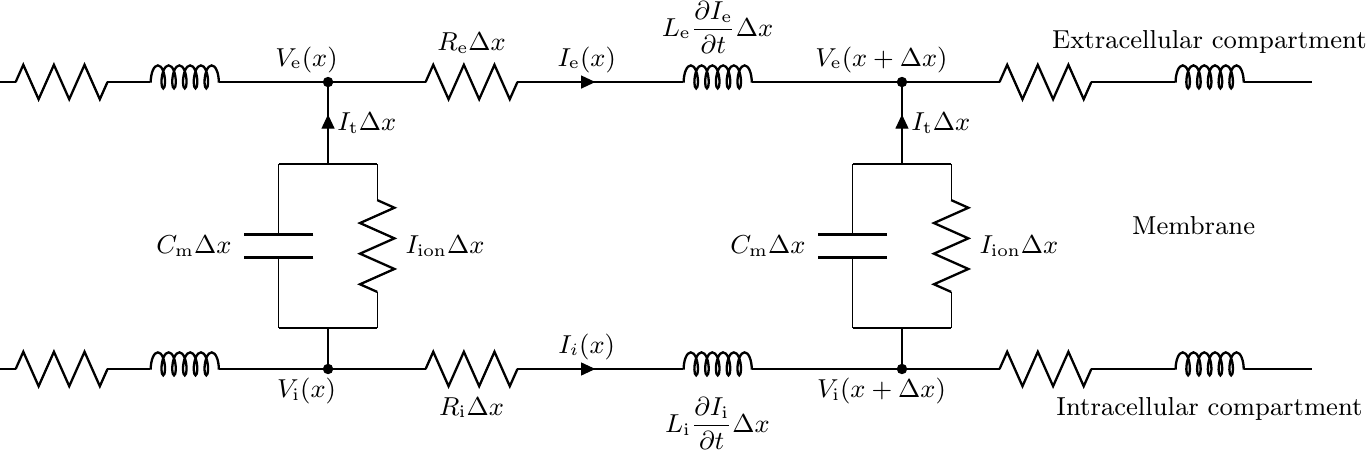}
\protect\caption{Schematic diagram of a discretized cable including inductance effects,
with isopotential circuit elements of length $\Delta x$. The inductances
in the cable give a Cattaneo-type model of the fluxes of the form
(\ref{eq:cattaneofluxe}) and (\ref{eq:cattaneofluxi}).}
\label{Fig:circuit}
\end{figure*}

\section{Phenomenological derivation of the hyperbolic bidomain equations\label{sec:bidomain}}

In the 1970's, \citet{Tung1978} formulated a bidomain model of the propagation of the action potential in cardiac muscle.
This tissue-scale model considers the myocardium to be composed of continuous intracellular and extracellular compartments, coupled via a continuous cellular membrane.
The bidomain equations can be derived from a model that accounts for the tissue microarchitecture \cite{Keener2009,NeuKrassowska1993,hand2009deriving}, but the bidomain model is a fundamentally homogenized description of excitation propagation that neglects the details of this microarchitecture.
Instead, the bidomain equations describe the dynamics of a local average of the voltages in the intracellular and extracellular compartments over a control volume.
One of the assumptions required by the homogenization procedure is that the control volume is large compared to the
scale of the cellular microarchitecture but small compared to any other important spatial scale of the system, such
as the width of the action potential wavefront.
Although the validity of this model has been questioned, for example by \citet{bueno2014fractional}, this approach has been extremely successful, and at present, most simulations of cardiac electrophysiology use such models.
For a detailed review of the bidomain model and other models of electrophysiology, we refer to \citet{griffith2013electrophysiology}.
Here, we assume that the homogenization assumptions hold, and we derive the hyperbolic bidomain model phenomenologically, starting from charge conservation in quasistatic conditions.

In the hyperbolic bidomain model, as in the standard bidomain model, the extracellular
and intracellular compartments are characterized by the anisotropic conductivity
tensors, respectively $\boldsymbol{\sigma}_{\text{e}}$ and $\boldsymbol{\sigma}_{\text{i}}$.
Introducing a local orthonormal basis, $\left\{
\boldsymbol{f}_{\text{0}},\boldsymbol{s}_{\text{0}},\boldsymbol{n}_{\text{0}}\right\} $, we assume the conductivity
tensors can be written as
$\boldsymbol{\sigma}_{j}=\sigma_{j}^{\text{f}}\boldsymbol{f}_{\text{0}}\otimes\boldsymbol{f}_{\text{0}}+\sigma_{j}^{\text{s}}\boldsymbol{s}_{\text{0}}\otimes\boldsymbol{s}_{\text{0}}+\sigma_{j}^{\text{n}}\boldsymbol{n}_{\text{0}}\otimes\boldsymbol{n}_{\text{0}},$
for $j=\mathrm{e},\,\mathrm{i}$. Typical values of the conductivity coefficients range from $10^{-2}$ mS/cm in the cross-fiber direction to 1 mS/cm in the fiber direction (for example, \citet{ColliFranzone2006} use $\sigma^{\text{f}}=1.2$
mS/cm, $\sigma^{\text{s}}=0.346$ mS/cm, and $\sigma^{\text{n}}=0.0435$ mS/cm). 
We model the current density fluxes using a Cattaneo-type equation, so that
\begin{eqnarray}
\tau_{\text{e}}\dfrac{\partial\boldsymbol{J}_{\text{e}}}{\partial t}+\boldsymbol{J}_{\text{e}} & = &
-\boldsymbol{\sigma}_{\text{e}}\nabla V_{\text{e}},\label{eq:cattaneofluxe}\\
\tau_{\text{i}}\dfrac{\partial\boldsymbol{J}_{\text{i}}}{\partial t}+\boldsymbol{J}_{\text{i}} & = &
-\boldsymbol{\sigma}_{\text{i}}\nabla V_{\text{i}},\label{eq:cattaneofluxi}
\end{eqnarray}
where $V_{\text{i}}$ ($V_{\text{e}}$), $\boldsymbol{J}_{\text{i}}$ ($\boldsymbol{J}_{\text{e}}$),
and $\tau_{\text{i}}$ ($\tau_{\text{e}}$) are the intracellular (extracellular)
potential, the intracellular (extracellular) flux, and the intracellular
(extracellular) relaxation time, respectively. As shown by \citet{Jou1996},
relations (\ref{eq:cattaneofluxe}) and (\ref{eq:cattaneofluxi})
can be also derived from the simplest model of ionic conduction in
a dilute system. The derivation of these evolution laws for the fluxes
from the generalized Gibbs equation \cite{Jou1996,Jou2002} shows
that they are consistent with the theory of extended irreversible
thermodynamics. Alternatively, equations (\ref{eq:cattaneofluxe})
and (\ref{eq:cattaneofluxi}) can be interpreted as arising from
a circuit model that includes inductances, such as the one depicted
in Fig.~\ref{Fig:circuit}. 
The derivation of the cable equation for
the circuit model shown in Fig.~\ref{Fig:circuit} is performed
in Appendix \ref{Appendix:Cable-equation}. In particular, we show in equation (\ref{eq:relaxtime})
the relationship between inductance and the relaxation time. 

To derive the higher-dimensional model equations, we begin by taking the divergence of the fluxes (\ref{eq:cattaneofluxe})
and (\ref{eq:cattaneofluxi}), so that
\begin{eqnarray}
\tau_{\text{e}}\dfrac{\partial}{\partial t}\nabla\cdot\boldsymbol{J}_{\text{e}}+\nabla\cdot\boldsymbol{J}_{\text{e}} & = & -\nabla\cdot\boldsymbol{\sigma}_{\text{e}}\nabla V_{\text{e}},\label{eq:divcattaneofluxe}\\
\tau_{\text{i}}\dfrac{\partial}{\partial t}\nabla\cdot\boldsymbol{J}_{\text{i}}+\nabla\cdot\boldsymbol{J}_{\text{i}} & = & -\nabla\cdot\boldsymbol{\sigma}_{\text{i}}\nabla V_{\text{i}}.\label{eq:divcattaneofluxi}
\end{eqnarray}
As in similar derivations of the bidomain model, we impose a quasistatic
form of charge conservation \cite{Sachse2004}, yielding
\begin{equation}
\nabla\cdot\left(\boldsymbol{J}_{\text{i}}+\boldsymbol{J}_{\text{e}}\right)=0.
\end{equation}
Additionally, the current leaving each compartment needs to enter the other, so that 
\begin{equation}
-\nabla\cdot\boldsymbol{J}_{\text{i}}=I_{\text{t}}=\nabla\cdot\boldsymbol{J}_{\text{e}},\label{eq:Itransmembrane}
\end{equation}
where $I_{\text{t}}=\chi\left(C_{\text{m}}\dfrac{\partial V}{\partial t}+I_{\text{ion}}\right)$
is the usual transmembrane current density, with $\chi$ the membrane area per unit volume of tissue, $C_{\text{m}}$ the membrane capacitance, and $I_{\text{ion}}$ the transmembrane ionic current. Using equation (\ref{eq:Itransmembrane})
in (\ref{eq:divcattaneofluxe}) and (\ref{eq:divcattaneofluxi}),
we obtain the system of equations,
\begin{eqnarray}
\tau_{\text{e}}\dfrac{\partial I_{\text{t}}}{\partial t}+I_{\text{t}} & = & -\nabla\cdot\boldsymbol{\sigma}_{\text{e}}\nabla V_{\text{e}},\label{eq:H0bidomaine}\\
-\tau_{\text{i}}\dfrac{\partial I_{\text{t}}}{\partial t}-I_{\text{t}} & = & -\nabla\cdot\boldsymbol{\sigma}_{\text{i}}\nabla V_{\text{i}}.\label{eq:H0bidomaini}
\end{eqnarray}
Defining the transmembrane potential $V=V_{\text{i}}-V_{\text{e}}$ to eliminate
$V_{\text{i}}$ from the equations yields
\begin{eqnarray}
\tau_{\text{e}}\dfrac{\partial I_{\text{t}}}{\partial t}+I_{\text{t}} & = & -\nabla\cdot\boldsymbol{\sigma}_{\text{e}}\nabla V_{\text{e}},\label{eq:bidomaingenformVe}\\
\tau_{\text{i}}\dfrac{\partial I_{\text{t}}}{\partial t}+I_{\text{t}} & = & \nabla\cdot\boldsymbol{\sigma}_{\text{i}}\nabla V+\nabla\cdot\boldsymbol{\sigma}_{\text{i}}\nabla V_{\text{e}}.\label{eq:bidomaingenformV}
\end{eqnarray}
Expanding the transmembrane currents, we finally obtain the hyperbolic
bidomain model,
\begin{eqnarray}
\tau_{\text{e}}C_{\text{m}}\dfrac{\partial^{2}V}{\partial t^{2}}+C_{\text{m}}\dfrac{\partial V}{\partial t}+\nabla\cdot\boldsymbol{D}_{\text{e}}\nabla V_{\text{e}}=\qquad\qquad\qquad\qquad\nonumber \\
-I_{\text{ion}}-\tau_{\text{e}}\dfrac{\partial I_{\text{ion}}}{\partial t},\qquad\label{eq:HHbidomainH1}\\
\tau_{\text{i}}C_{\text{m}}\dfrac{\partial^{2}V}{\partial t^{2}}+C_{\text{m}}\dfrac{\partial V}{\partial t}\qquad\qquad\qquad\qquad\qquad\qquad\qquad\nonumber \\
-\nabla\cdot\boldsymbol{D}_{\text{i}}\nabla V-\nabla\cdot\boldsymbol{D}_{\text{i}}\nabla V_{\text{e}}=\qquad\qquad\nonumber \\
-I_{\text{ion}}-\tau_{\text{i}}\dfrac{\partial I_{\text{ion}}}{\partial t},\qquad\label{eq:HHbidomainH2}
\end{eqnarray}
where we set $\boldsymbol{D}_{\text{i}}=\boldsymbol{\sigma}_{\text{i}}/\chi$ and $\boldsymbol{D}_{\text{e}}=\boldsymbol{\sigma}_{\text{e}}/\chi$. Alternatively,
the model equations can be written as
\begin{eqnarray}
\tau_{\text{i}}C_{\text{m}}\dfrac{\partial^{2}V}{\partial t^{2}}+C_{\text{m}}\dfrac{\partial V}{\partial t}\qquad\qquad\qquad\qquad\qquad\qquad\qquad\nonumber \\
-\nabla\cdot\boldsymbol{D}_{\text{i}}\nabla V-\nabla\cdot\boldsymbol{D}_{\text{i}}\nabla V_{\text{e}}=\qquad\qquad\nonumber \\
-I_{\text{ion}}-\tau_{\text{i}}\dfrac{\partial I_{\text{ion}}}{\partial t},\qquad\label{eq:HH2bidomainH1}\\
\left(\tau_{\text{e}}-\tau_{\text{i}}\right)C_{\text{m}}\dfrac{\partial^{2}V}{\partial t^{2}}\qquad\qquad\qquad\qquad\qquad\qquad\qquad\qquad\nonumber \\
+\nabla\cdot\boldsymbol{D}_{\text{i}}\nabla V+\nabla\cdot\left(\boldsymbol{D}_{\text{e}}+\boldsymbol{D}_{\text{i}}\right)\nabla V_{\text{e}}=\qquad\qquad\nonumber \\
\left(\tau_{\text{i}}-\tau_{\text{e}}\right)\dfrac{\partial I_{\text{ion}}}{\partial t}.\qquad\label{eq:HH2bidomainH2}
\end{eqnarray}
Notice that if $\tau_{\text{i}}=\tau_{\text{e}}=\tau$, then these equations reduce
to a hyperbolic-elliptic system that is similar to the parabolic-elliptic form of the standard bidomain model,
\begin{eqnarray}
\tau C_{\text{m}}\dfrac{\partial^{2}V}{\partial t^{2}}+C_{\text{m}}\dfrac{\partial V}{\partial t}\qquad\qquad\qquad\nonumber \\
-\nabla\cdot\boldsymbol{D}_{\text{i}}\nabla V-\nabla\cdot\boldsymbol{D}_{\text{i}}\nabla V_{\text{e}} & = & -I_{\text{ion}}-\tau\dfrac{\partial I_{\text{ion}}}{\partial t},\quad\label{eq:HEbidomainH}\\
\nabla\cdot\boldsymbol{D}_{\text{i}}\nabla V+\nabla\cdot\left(\boldsymbol{D}_{\text{e}}+\boldsymbol{D}_{\text{i}}\right)\nabla V_{\text{e}} & = & 0.\label{eq:HEbidomainE}
\end{eqnarray}
If we further take $\tau=0$, we retrieve the usual bidomain model
in its parabolic-elliptic form,
\begin{eqnarray}
C_{\text{m}}\dfrac{\partial V}{\partial t}-\nabla\cdot\boldsymbol{D}_{\text{i}}\nabla V-\nabla\cdot\boldsymbol{D}_{\text{i}}\nabla V_{\text{e}} & = &- I_{\text{ion}},\label{eq:Pbidomain}\\
\nabla\cdot\boldsymbol{D}_{\text{i}}\nabla V+\nabla\cdot\left(\boldsymbol{D}_{\text{e}}+\boldsymbol{D}_{\text{i}}\right)\nabla V_{\text{e}} & = & 0.\label{eq:Ebidomain}
\end{eqnarray}

\section{Reduction to the hyperbolic monodomain \label{sec:monodomain}}

The hyperbolic bidomain model can be simplified by assuming the extracellular and intracellular compartments have the same anisotropy ratios, so that $\boldsymbol{D}=\boldsymbol{D}_{\text{e}}=\lambda\boldsymbol{D}_{\text{i}}$.
If we make this assumption, we obtain from (\ref{eq:HH2bidomainH2})
\begin{eqnarray}
-\nabla\cdot\boldsymbol{D}\nabla V_{\text{e}}=\dfrac{1}{\lambda+1}\nabla\cdot\boldsymbol{D}\nabla V\qquad\qquad\qquad\qquad\qquad\nonumber \\
+\dfrac{\lambda}{\lambda+1}\left(\tau_{\text{e}}-\tau_{\text{i}}\right)C_{\text{m}}\dfrac{\partial^{2}V}{\partial t^{2}}-\dfrac{\lambda}{\lambda+1}\left(\tau_{\text{i}}-\tau_{\text{e}}\right)\dfrac{\partial I_{\text{ion}}}{\partial t},\qquad\label{eq:Hmonodomain1}
\end{eqnarray}
Substituting in (\ref{eq:HH2bidomainH1}), we find
\begin{eqnarray}
\left[\tau_{\text{i}}+\dfrac{\lambda}{\lambda+1}\left(\tau_{\text{e}}-\tau_{\text{i}}\right)\right]C_{\text{m}}\dfrac{\partial^{2}V}{\partial t^{2}}\qquad\qquad\qquad\qquad\nonumber \\
+C_{\text{m}}\dfrac{\partial V}{\partial t}-\dfrac{\lambda}{\lambda+1}\nabla\cdot\boldsymbol{D}\nabla V=-I_{\text{ion}}\qquad\qquad\nonumber \\
-\left[\tau_{\text{i}}+\dfrac{\lambda}{\lambda+1}\left(\tau_{\text{e}}-\tau_{\text{i}}\right)\right]\dfrac{\partial I_{\text{ion}}}{\partial t}.\quad\label{eq:Hmonodomain2}
\end{eqnarray}
Defining $\tau=\tau_{\text{i}}+\lambda\left(\tau_{\text{e}}-\tau_{\text{i}}\right)/\left(\lambda+1\right)$,
we obtain the hyperbolic monodomain model,
\begin{equation}
\tau C_{\text{m}}\dfrac{\partial^{2}V}{\partial t^{2}}+C_{\text{m}}\dfrac{\partial V}{\partial t}-\nabla\cdot\boldsymbol{D}\nabla V=-I_{\text{ion}}-\tau\dfrac{\partial I_{\text{ion}}}{\partial t},\label{eq:Hmonodomain}
\end{equation}
where here we have absorbed the term $\lambda/\left(\lambda+1\right)$ into $\boldsymbol{D}$.
The relaxation time $\tau$ of the monodomain model is always positive, and
it is zero only if both $\tau_{\text{i}}$ and $\tau_{\text{e}}$ are zero. In fact,
if $\tau_{\text{e}}=0$, then $\tau=\tau_{\text{i}}/\left(\lambda+1\right)$, and if
$\tau_{\text{i}}=0$, then $\tau=\lambda\tau_{\text{e}}/\left(\lambda+1\right)$.
However, if $\tau_{\text{i}}=\tau_{\text{e}}$, then $\tau=\tau_{\text{i}}=\tau_{\text{e}}$. 

Introducing a new variable $Q = \dfrac{\partial V}{\partial t}$, we transform the hyperbolic monodomain
equations into the first-order system, 
\begin{eqnarray}
\dfrac{\partial V}{\partial t} & = & Q,\label{eq:Vequation}\\
\tau C_{\text{m}}\dfrac{\partial Q}{\partial t}+C_{\text{m}}Q-\nabla\cdot\boldsymbol{D}\nabla V & = & -I_{\text{ion}}-\tau\dfrac{\partial I_{\text{ion}}}{\partial t}.\qquad\label{eq:Qequation}
\end{eqnarray}
Equations (\ref{eq:Vequation})--(\ref{eq:Qequation}) are usually
supplemented with insulation boundary conditions, such that 
\begin{equation}
\nabla V\cdot\boldsymbol{N}=0\label{eq:BC}
\end{equation}
 on the boundary of the domain, where the vector $\boldsymbol{N}$
is the normal to the boundary.

Following \citet{stan2014finite}, we say that a solution of the hyperbolic monodomain equations has a finite propagation speed if, given compactly supported initial conditions for $V$ at time $t=0$, $V(\cdot,t)$ is also compactly supported for any $t > 0$. The compact support is taken with respect to the resting potential $V_0$.
By contrast, a solution has infinite propagation speed if the initial data are compactly supported, but for any $t > 0$ and any $R > 0$, the set 
$ M_{R,t} = \{\boldsymbol{x}  : ||\boldsymbol{x}||_2 \geq R, V(\boldsymbol{x}, t) > V_0 \} $ has positive measure.
For any  relaxation time $\tau\neq0$, (\ref{eq:Hmonodomain}) is  hyperbolic, and it thereby has the property of finite propagation speed. Setting $\tau = 0$, the equations become parabolic and have solutions with infinite propagation speeds. Specifically, in the standard parabolic model, local perturbations in $V$, even those that do not generate a propagating front, will travel at infinite speed.

\section{Ionic models\label{sec:Ionicmodels}}

\begin{table}[t]
\centering
\begin{tabular*}{1\columnwidth}{@{\extracolsep{\fill}}|ccccccccc|}
\hline 
$C_{\text{m}}$ & $\sigma^{\text{f}}$ & $\sigma^{\text{s}}=\sigma^{\text{n}}$ & $\chi$ & $k$ & $b$ & $\mu_{1}$ & $\mu_{2}$ & $\varepsilon$\tabularnewline
\hline 
1 & 1 & 0.125 & 1 & 8 & 0.1 & 0.12 & 0.3 & 0.01\tabularnewline
\hline 
\end{tabular*}
\caption{Parameters for the Aliev-Panfilov ionic model (\ref{eq:AP}).}
\label{Tab:NP}
\end{table}

\begin{table*}[t]
\centering
\begin{tabular*}{1\textwidth}{@{\extracolsep{\fill}}|c|ccccccccccccccccc|}
\hline 
 & $C_{\text{m}}$ & $\sigma^{\text{f}}$ & $\sigma^{\text{s}}=\sigma^{\text{n}}$ & $\chi$ & $\tau_{v}^{+}$ & $\tau_{v1}^{-}$ & $\tau_{v2}^{-}$ & $\tau_{w}^{+}$ & $\tau_{w}^{-}$ & $\tau_\text{d}$ & $\tau_{0}$ & $\tau_\text{r}$ & $\tau_\text{si}$ & $k$ & $V_\text{c}^\text{si}$ & $V_\text{c}$ & $V_v$\tabularnewline
\hline 
parameter set 3 & 1 & 0.1 & 0.0125 & 1 & 3.33 & 19.6 & 1250 & 870 & 41 & 0.25 & 12.5 & 33.33 & 29.0 & 10 & 0.85 & 0.13 & 0.04\tabularnewline
parameter set 4 & 1 & 0.1 & 0.0125 & 1 & 3.33 & 15.6 & 5 & 350 & 80 & 0.407 & 9.0 & 34.0 & 26.5 & 15 & 0.45 & 0.15 & 0.04\tabularnewline
parameter set 5 & 1 & 0.1 & 0.0125 & 1 & 3.33 & 12.0 & 2 & 1000 & 100 & 0.362 & 5.0 & 33.33 & 29.0 & 15 & 0.7 & 0.13 & 0.04\tabularnewline
parameter set 6 & 1 & 0.1 & 0.0125 & 1 & 3.33 & 9.0 & 8 & 250 & 60 & 0.395 & 9.0 & 33.33 & 29.0 & 15 & 0.5 & 0.13 & 0.04\tabularnewline
\hline 
\end{tabular*}
\caption{Sets of parameters for the Fenton-Karma ionic model (\ref{eq:FKv})--(\ref{eq:Isi}) as stated by \citet{Fenton2002}.}
\label{Tab:FK}
\end{table*}

Transmembrane ionic fluxes through ion channels, pumps, and exchangers are responsible for the cardiac action potential.
The action potential is initiated by a fast inward sodium current that depolarizes the cellular membrane.
After depolarization phase, slow inward currents (primarily calcium currents) and slow outward currents (primarily potassium currents) approximately balance each other, prolonging the action potential and creating a plateau phase.
Ultimately, the slow outward currents bring the transmembrane potential difference back to its resting value of approximately $-80~\text{mV}$.
The bidomain and monodomain models must be completed by specifying the form of $I_{\text{ion}}$, which accounts for these transmembrane currents. 
The states of the transmembrane ion channels are described by a collection of variables, $\boldsymbol{w}$, associated with the ionic model, so that $I_{\text{ion}} = I_{\text{ion}}(V,\boldsymbol{w})$.
Typically, the state variable dynamics are determined by a spatially decoupled system of nonlinear ordinary differential equations,
\begin{equation}
\dfrac{\partial\boldsymbol{w}}{\partial t}=\boldsymbol{g}\left(V,\boldsymbol{w}\right),\label{eq:ionicmodel}
\end{equation}
where the form of $\boldsymbol{g}$ depends on the details of the particular ionic model.
We consider the simplified piecewise-linear model of \citet{McKean1970}, the two-variable model of \citet{Aliev1996},
the three-variable model of \citet{Fenton1998}, and the biophysically
detailed models of \citet{TenTusscher2006} for the ventricles (20
variables) and of \citet{Grandi2011} for the atria (57 variables). 

In the two-variable model of \citet{Aliev1996}, $\boldsymbol{w}=\left\{ r\right\} $, with
\begin{equation}
\dfrac{\partial r}{\partial t}=\left(\varepsilon+\dfrac{\mu_{1}r}{\mu_{2}+V}\right)\left(-r-kV\left(V-b-1\right)\right),\label{eq:AP}
\end{equation}
and the ionic current is $I_{\text{ion}}=kV\left(V-\alpha\right)\left(V-1\right)+rV$.
The parameters used in this case are shown in Table \ref{Tab:NP}.
In the three-variable model of \citet{Fenton1998}, $\boldsymbol{w}=\left\{ v,w\right\} $,
with
\begin{eqnarray}
\dfrac{\partial v}{\partial t} & = & \dfrac{1}{\tau_{v}^{-}\left(V\right)}\left(1-p\right)\left(1-v\right)-\dfrac{1}{\tau_{v}^{_{+}}}pv,\label{eq:FKv}\\
\dfrac{\partial w}{\partial t} & = & \dfrac{1}{\tau_{w}^{-}}\left(1-p\right)\left(1-w\right)-\dfrac{1}{\tau_{w}^{_{+}}}pw,\label{eq:FKw}
\end{eqnarray}
where $\tau_{v}^{-}\left(V\right)=\left(1-q\right)\tau_{v1}^{-}+q\tau_{v2}^{-}$, $p=H\left(V-V_\text{c}\right)$, $q=H\left(V-V_{v}\right)$, and $H\left(\cdot\right)$ is the Heaviside function.
The total
ionic current is the sum of three currents, $I_{\text{ion}}=-I_\text{fi}\left(V,v\right)-I_\text{so}\left(V\right)-I_\text{si}\left(V,w\right)$,
with
\begin{eqnarray}
I_\text{fi}\left(V,v\right) & = & -\dfrac{1}{\tau_\text{d}}vp\left(V-V_\text{c}\right)\left(1-V\right),\label{eq:Ifi}\\
I_\text{so}\left(V\right) & = & \dfrac{1}{\tau_{0}}V\left(1-p\right)+\dfrac{1}{\tau_\text{r}}p,\label{eq:Iso}\\
I_\text{si}\left(V,w\right) & = & -\dfrac{1}{2\tau_\text{si}}w\left(1+\tanh\left(k\left(V-V_\text{c}^\text{si}\right)\right)\right).\label{eq:Isi}
\end{eqnarray}
The sets of parameters used in this paper for this model are shown in Table \ref{Tab:FK}. 
We refer to the original papers \cite{TenTusscher2006, Grandi2011} for the statements of the equations and parameters
of the biophysically detailed ionic models.

\begin{figure*}
\begin{centering}
\includegraphics{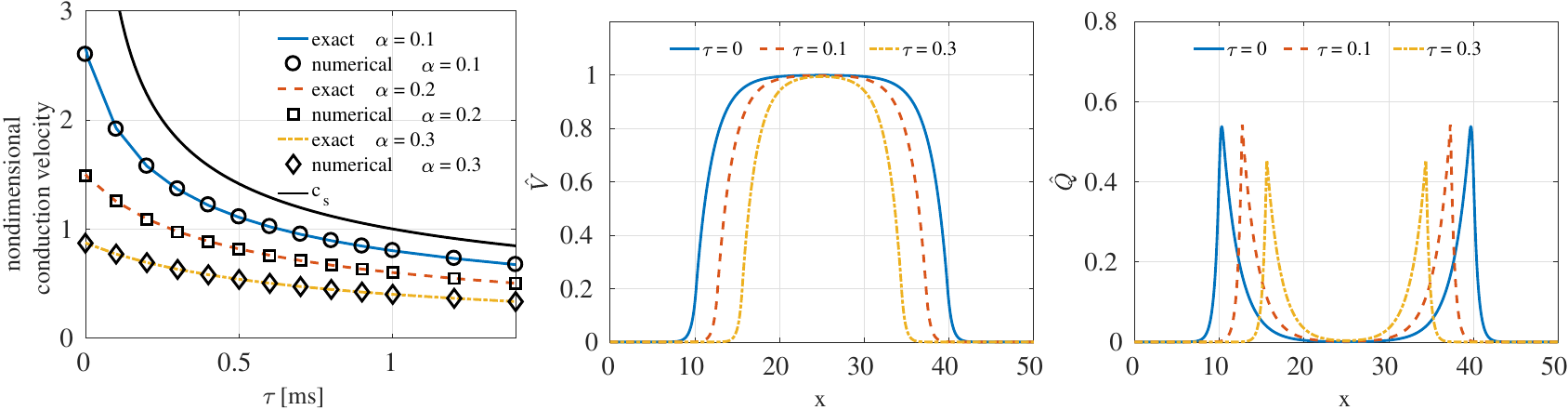}
\par\end{centering}

\protect\caption{Conduction velocities for the piecewise-linear reaction model (\ref{eq:mckeanpiecewiselinear}).
Left) comparison between the numerical and the exact wave speed (\ref{eq:mckeanspeed}).
The black curve represent the characteristic propagation speed of the system,
i.e.~the limit speed at which local perturbations can travel. If the relaxation
time $\tau=0$, then local perturbations can propagate at infinite
speed. Center) Propagating front for different relaxation times at
$t=10$. Right) time derivative of the nondimensional potential for
different relaxation times at $t=10$. All simulations were run using a
mesh size $h=0.03125$ and a time step size $\Delta t=0.003653.$}

\label{Fig:bistable}
\end{figure*}

\section{Numerical Results\label{sec:NumericalResults}}

The hyperbolic monodomain and bidomain models
are discretized using a low-order finite element scheme described in Appendix \ref{sec:FEM} that is implemented in
the open-source parallel C++ code BeatIt (available at \url{http://github.com/rossisimone/beatit}),
which is based on the libMesh finite element library\citep{libMeshPaper} and relies on
linear solvers provided by PETSc \citep{petsc-efficient,petsc-user-ref,petsc-web-page}.
All the code used for the following tests is contained in the
online repository, and all tests can be replicated directly from those codes. The only exception is the atrial fibrillation test, which 
uses a patient-specific mesh that is not contained in the repository.
One-~and two-dimensional simulations were run using a Linux workstation with two
Intel Xeon E5-2650 v3 processors (up to 40 threads) and 32 GB
of memory. Three-dimensional simulations were run on the KillDevil Linux cluster
at the University of North Carolina at Chapel Hill. We used Matlab\citep{MATLAB:2010}
to visualize the one-dimensional results and Paraview\citep{ahrens200536} for
the two-~and three-dimensional simulations.

\subsection{Comparison with an exact solution\label{sub:resultsexact}}

We start by considering a simple piecewise-linear bistable model for the ionic currents,
with
\begin{eqnarray}
I_{\text{ion}}\left(V\right) & = & kV\nonumber \\
 & - & k\left[V_{2}H(V-V_{1})+V_{0}\left(1-H(V-V_{1})\right)\right],\quad\label{eq:dimmckean}
\end{eqnarray}
where $V_{0}$ is the resting potential, $V_{1}$ is the threshold
potential, $V_{2}$ is the depolarization potential, and $H\left(\cdot\right)$ is the Heaviside function. This model (\ref{eq:dimmckean})
reduces to the piecewise-linear model of \citet{McKean1970} after
a simple dimensional analysis, after which the nondimensional ionic currents
take the form
\begin{equation}
\hat{I}_{\text{ion}}\left(\hat{V}\right)=\hat{V}-H\left(\hat{V}-\alpha\right).\label{eq:mckeanpiecewiselinear}
\end{equation}
The nondimensional
potential $\hat{V}$ is related to the dimensional potential by
$V=\left(V_{2}-V_{0}\right)\hat{V}+V_{0}$.  In this model,
$\alpha=\left(V_{1}-V_{0}\right)/\left(V_{2}-V_{0}\right)$ describes the excitability of the tissue.
As we show in Appendix \ref{Appendix:exactsolution}, using model
(\ref{eq:mckeanpiecewiselinear}), it is possible to find the analytic
solution of a propagating front for the nondimensional hyperbolic monodomain problem. In
particular, we find that the front propagation speed in an unbounded domain is
\begin{equation}
c=\dfrac{\left(1-2\alpha\right)}{\sqrt{\mu+\left(\alpha-\alpha^{2}\right)\left(\mu-1\right)^{2}}}<\sqrt{\dfrac{1}{\mu}}=c_\text{s}.\label{eq:mckeanspeed}
\end{equation}
The speed $c_\text{s}=\sqrt{1/\mu}=\sqrt{C_{\text{m}}/\tau k}$ represents the
maximum speed at which a perturbation can travel in the system. When
the relaxation time $\tau$ goes to zero, the local perturbations can
travel at infinite speed. The parameter $\mu = \tau k / C_{\text{m}}$ is a nondimensional
number representing the ratio between the relaxation time and the characteristic time of the reactions and characterizing the effects of the inductances in the system (see Appendix \ref{Appendix:exactsolution}).
The corresponding dimensional speed
is
\begin{equation}
v=\sqrt{\dfrac{\sigma k}{\chi C_{\text{m}}^{2}}}\dfrac{\left(1-2\alpha\right)}{\sqrt{\mu+\left(\alpha-\alpha^{2}\right)\left(\mu-1\right)^{2}}}<\sqrt{\dfrac{\sigma}{\chi\tau C_{\text{m}}}}.\label{eq:mckeanspeedphysical}
\end{equation}

As a verification test, we consider the spatial interval $\Omega=\left[0,50\right]$. An initial
stimulus of amplitude $1$ is applied at $x\in\left[24.5,25.5\right]$
for the interval $t\in\left[0.03,1.03\right]$. The system of equations is solved
using a mesh size $h=0.03125$ and a time step size $\Delta t=0.003653.$
We show the nondimensional solutions $\hat{V}$ and $\hat{Q}$ at
$t=10$ in Fig.~\ref{Fig:bistable}(center and right). We register
the activation time at a particular spatial location whenever $\hat{V}$ crosses a threshold of $0.9$. The conduction velocities
are measured by picking the distance between two points $x_{1}$ and
$x_{2}$ $\in\Omega$ and dividing it by the time interval between
the activation times at these two locations. Because (\ref{eq:mckeanspeed})
is obtained by assuming that $\Omega$ is the entire real line, we need to measure
conduction velocities far from the boundaries. For this reason, we
choose $x_{1}=30$ and $x_{2}=32$. The comparison between the exact
conduction velocities defined by equation (\ref{eq:mckeanspeed})
and those obtained by the simulations is shown in Fig.~\ref{Fig:bistable}(left).

We also examine the effect of the relaxation time on the conduction
velocities for three values of the excitability parameter, $\alpha\in\left\{ 0.1,\,0.2,\,0.3\right\} $.
With this simple ionic model, the larger the relaxation time, the slower the wave speed. The
limit speed $c_\text{s}$ at which local perturbations can travel is also
shown in Fig.~\ref{Fig:bistable}(left). Thus, in this piecewise-linear model,
the effect of inductances is to slow down the propagation of the front.
By contrast, the following tests will show that with fully nonlinear models, inductances can also enhance propagation.

\begin{figure*}
\begin{centering}
\includegraphics{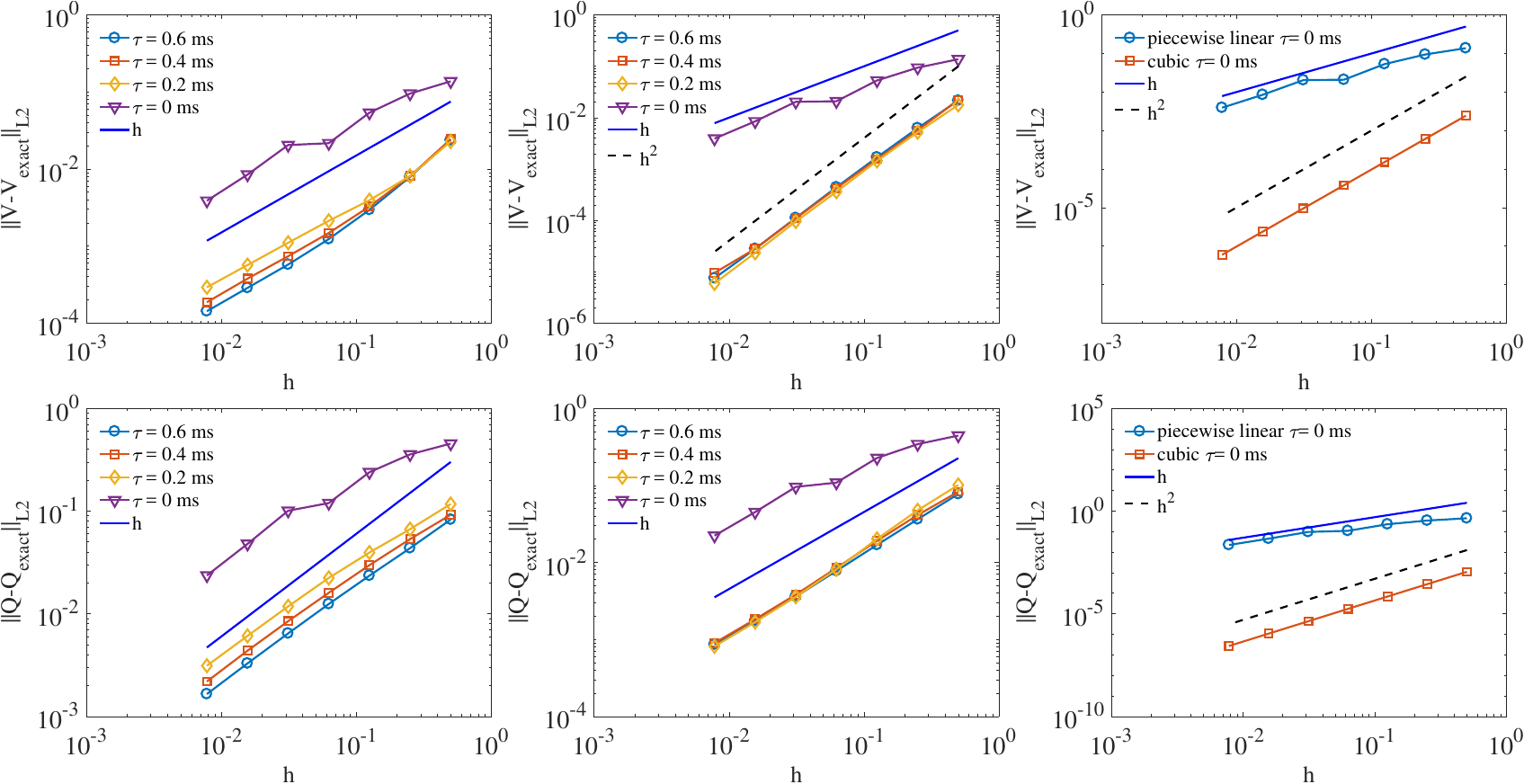}

\par\end{centering}

\protect\caption{Convergence study for the piecewise-linear 
reaction model (\ref{eq:mckeanpiecewiselinear}) for the variables $V$ and $Q$.
Left) comparison of the rate of convergence for the hyperbolic and parabolic
monodomain models using a first-order implicit-explicit (IMEX) Runge-Kutta (RK) time integrator. Center) comparison of the rate of convergence for the hyperbolic and parabolic
monodomain models
using a second-order IMEX-RK time integrator. Right) comparison of the rate of
convergence between the piecewise-linear reaction model and a cubic reaction
 model for the parabolic monodomain model using a second-order IMEX-RK time
integrator. }

\label{Fig:bistable_convergence}
\end{figure*}

Using the analytic solution derived in Appendix \ref{Appendix:exactsolution}, we also perform a
space-time convergence study. We consider the spatial interval $\Omega=[-25,25]$ and set the initial conditions
according to
\eqref{eq:continuituconditions} and \eqref{eq:Derivative}, assuming the front at time $t=0$ is at $x=0$ and $\alpha =
0.1$.
The system is discretized in time using implicit-explicit (IMEX) Runge-Kutta (RK) time integrators\cite{ascher1997implicit,boscarino2016high}. In
particular, we use the forward/backward Euler schemes for the first-order time integrator and the Heun/Crank-Nicholson schemes for the second-order time integrators. For more details on the time integrator, see Appendix \ref{sec:FEM}. To ensure the validity of
the solutions \eqref{eq:continuituconditions} and \eqref{eq:Derivative} in the considered bounded domain, we run the simulation from $t=0$ to $t=1$, so that the front is still far from the boundaries. 
The time step size $\Delta t$ is taken to be linearly proportional to the mesh size $h$, with a value of $\Delta t = 0.05$ in the coarsest cases.
The time step size was chosen to guarantee a large enough number of time iterations and while satisfying the CFL
condition (see Appendix \ref{sec:FEM}).
Fig.~\ref{Fig:bistable_convergence} shows the errors for the potential $V$ and its time derivative $Q$ using the first-order
(Fig.~\ref{Fig:bistable_convergence} left) and second-order (Fig.~\ref{Fig:bistable_convergence} center) time stepping schemes. Because the ionic currents and their
derivative in this case are not smooth functions, we do not expect to obtain full second-order convergence. 
 On the other hand, whereas the convergence rates for the parabolic monodomain model are always
first order in both $V$ and $Q$, the hyperbolic monodomain model with the second-order time stepping scheme converges quadratically in $V$ and linearly only in $Q$. 
The simulations were run with and without regularization of the Heaviside and Dirac-$\delta$ functions. The errors
reported in Fig. \ref{Fig:bistable_convergence}  correspond to the more accurate simulations.
We show in Fig.~\ref{Fig:bistable_convergence}(right) that
the slow convergence results from the non-smoothness of the ionic currents: we compare the convergence of the
piecewise-linear reaction model with the monodomain model with the cubic reaction term
\[\hat{I}_{\mathrm{ion}}(\hat{V})=\hat{V}(\hat{V}-1)(\hat{V}-\alpha).\]
The exact solution for the cubic reaction in the parabolic monodomain has been determined previously \cite{Keener2009,Pezzuto2016}. We are
not aware of the existence of an analytical solution for the hyperbolic monodomain with a cubic reaction term. We show
in Fig.~\ref{Fig:bistable_convergence} that the cubic model converges quadratically (using a second-order time stepping
scheme) whereas the piecewise-linear model converges linearly.
To conclude, at least in these tests, the discretization errors of the parabolic mondomain model are larger than the discretization errors in the hyperbolic monodomain model,
indicating a better accuracy for the hyperbolic model.

\begin{figure*}[htb]
\begin{centering}
\includegraphics{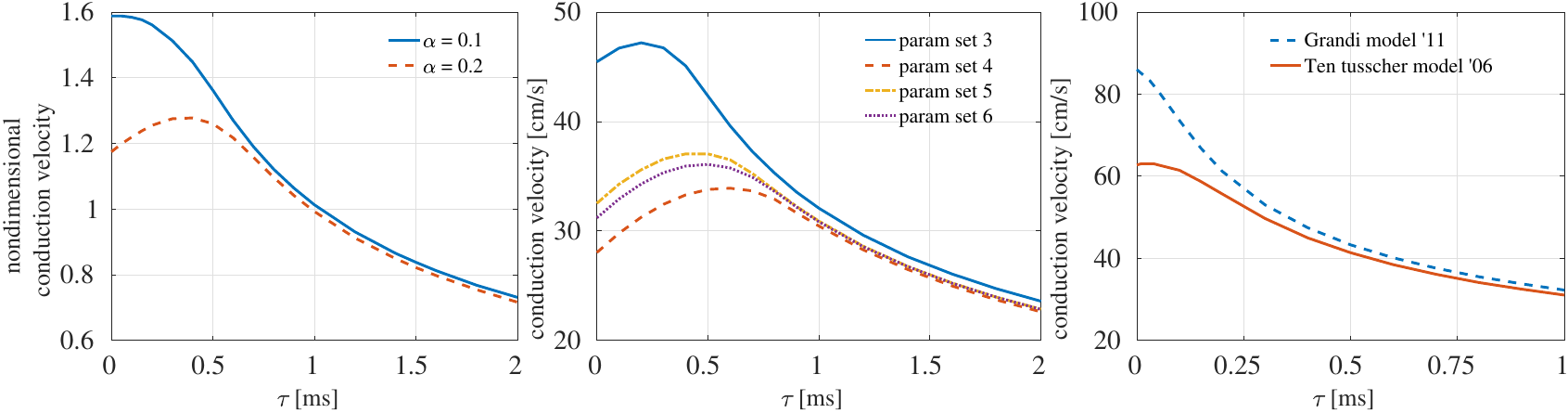}
\par\end{centering}

\protect\caption{Dependence of the conduction velocity on the relaxation time for different
ionic models: left) Aliev-Panfilov; center) Fenton-Karma; right) ventricular
Ten tusscher '06 model and atrial Grandi '11 model. An external stimulus
is applied at the center of a one-dimensional domain and the activation
times are recorded. For all models we used a time step size $\Delta t=0.003125$
ms and a mesh size $h=31.25\text{ }\mu$m over a domain of length
$5$ cm. Unlike the linear case, the nonlinearities of the time derivative
of the ionic currents increase the wavefront speed for small and moderate
values of the relaxation time $\tau$.}
\label{fig:test2}
\end{figure*}

\subsection{Effect of the relaxation time on the conduction velocity for different
ionic models\label{sub:resultsionicmodels}}

We now analyze how the inductance terms influence the conduction velocities
of cardiac action potentials. We consider four ionic models: the two-variable model of \citet{Aliev1996}, the three-variable model of
\citet{Fenton1998}, the ventricular M-cell model of \citet{TenTusscher2006}
(20 variables), and the atrial model of \citet{Grandi2011} (57 variables).
For the two-variable Aliev-Panfilov model, we use the parameters
of \citet{Nash2004}, as reported in Table \ref{Tab:NP}. We consider
two values for the excitability parameter, $\alpha=0.1$ and $\alpha=0.2$.
For the three-variable Fenton-Karma model, we consider parameters sets 3, 4,
5, and 6 given by \citet{Fenton2002} and reported
in Table \ref{Tab:FK}. For the biophysically detailed ionic model
of Ten tusscher et al.~and Grandi et al., our simulations used the C++ implementations of the models provided
by the authors, using $C_\text{m}=1$ $\mu\text{F}/\text{cm}^2$
for the membrane capacitance in both models. 

For all cases, the domain is the interval $\Omega=\left[0,5\right]$
cm. The simulations use a mesh size $h=31.25$ $\mu$m and a time step size
$\Delta t=0.003653$ ms. An initial stimulus is applied at the center
of the domain, $x\in\left[2.45,2.55\right]$ cm, during the time
interval $t\in\left[0.03,1.03\right]$ ms. The conduction velocities
are measured, as in the previous test case, by dividing the distance
between two points $x_{1}$ and $x_{2}$
by the time interval between the activation times at these two
selected locations. The computed conduction velocities of different ionic models
are shown in Fig.~\ref{fig:test2} for relaxation times in the range
$\left[0,1\right]$ ms. Contrary to the results for the piecewise-linear ionic model,
in this case, small and moderate values of the
relaxation time act to enhance propagation.

We also compare the computational time required by the parabolic and hyperbolic monodomain models in
two and three spatial dimensions. Table \ref{Table:cpu_time} shows the number of iterations taken by the linear solvers (conjugate gradient preconditioned by successive over-relaxation) along with the
relative computational times of the hyperbolic monodomain model (with respect to the parabolic model) for the \emph{reaction step} and
 the \emph{diffusion step}.
In the reaction step, we solve the ionic model, and we also evaluate the ionic currents and their time derivatives.
For this reason, we expect this step to take longer in the hyperbolic monodomain model. In fact, with simple ionic
models, the amount of work required to evaluate the ionic currents and their time derivatives is almost doubled. This is
natural because, for such simple models, the time derivative of the ionic currents do not use any quantities already computed in the evaluation of the ionic currents.
By contrast, for biophysically detailed models, most of the computation is needed for the evaluation of the ionic
currents, and many quantities can be reused for the evaluation of their time derivatives. This is reflected in Table
\ref{Table:cpu_time} by a small increase (less than 7\%)
in the computational time of the hyperbolic model for the Ten Tusscher '06 model. 
Unexpectedly, the solution of the linear system take fewer iterations in the hyperbolic monodomain model. This is
reflected by a speedup of about 20\% 
of the computational time in the diffusion step, where we assembled the right hand side, solve the linear system, and
update the variables $Q$ and $V$. 

\begin{table*}[t]
\begin{tabular*}{1\textwidth}{@{\extracolsep{\fill}}|c|c|c|c|c|c|}
\hline 
2D & \multicolumn{2}{c|}{Max CG Iterations} & Reaction Step Cost& Diffusion Step Cost& Overall\tabularnewline
\hline 
100$\times$100 & $\tau=0$  ms& $\tau=0$.4  ms & $\tau=0.4 \text{ ms}/\tau=0$ ms & $\tau=0.4 \text{ ms}/\tau=0$ ms &
$\tau=0.4 \text{ ms}/\tau=0$ ms\tabularnewline
\hline 
\hline 
AP & 8 & 3 & -54\% & 23\% & -3\%\tabularnewline
\hline 
FK & 22 & 15 & -62\% & 22\% & -2\%\tabularnewline
\hline 
TP06 & 11 & 5 & -7\% & 15\% & -3\%\tabularnewline
\hline 
\multicolumn{1}{c}{} & \multicolumn{1}{c}{} & \multicolumn{1}{c}{} & \multicolumn{1}{c}{} & \multicolumn{1}{c}{} & \multicolumn{1}{c}{}\tabularnewline
\hline 
3D & \multicolumn{2}{c|}{Max CG Iterations} & Reaction Step Cost& Diffusion Step Cost& Overall\tabularnewline
\hline 
10$\times$10$\times$10 & $\tau=0$ ms& $\tau=0$.4  ms& $\tau=0.4 \text{ ms}/\tau=0$ & $\tau=0.4 \text{ ms}/\tau=0$ ms &
$\tau=0.4 \text{ ms}/\tau=0$  ms\tabularnewline
\hline 
\hline 
AP & 43 & 16 & -73\% & 48\% & 30\%\tabularnewline
\hline 
FK & 43 & 16 & -49\% & 47\% & 25\%\tabularnewline
\hline 
TP06 & 43 & 16 & -6\% & 46\% & 15\%\tabularnewline
\hline 
\end{tabular*}
\caption{ Comparison of the compuational cost for the parabolic
($\tau = 0$ ms) and hyperbolic ($\tau = 0.4$ ms) monodomain models for two-~and three-dimensional problems in a serial computation. Negative values indicate relative slowdown of the hyperbolic model compared to the parabolic model, whereas
positive values represent relative speedup. In the \emph{reaction step}, we measure the time needed to solve the ionic model
and to evaluate the ionic currents and their time derivatives nodewise for 1000 time steps. In the \emph{diffusion step}, we
measure the time needed to form the right hand side, solve the linear system and update the solutions, for 1000 time steps. The hyperbolic monodomain model can be solved using fewer linear solver iterations. This is
reflected by the fact that the diffusion step is about 20\%
faster than in the parabolic monodomain model. On the other hand, the solution of the hyperbolic monodomain model
requires the additional evaluation of the derivative of the ionic currents. For this reason, the reaction step is
actually slower than in the parabolic monodomain model. Notice, however, the
difference in the computational cost between the hyperbolic and parabolic
reaction steps is smaller for the more complex ionic model. AP: Aliev-Panfilov ionic model; FK: Fenton-Karma ionic model; TP06: TenTusscher et al.~'06 ionic
model.}
\label{Table:cpu_time}
\end{table*}

\begin{figure*}[htb]
\centering
\includegraphics{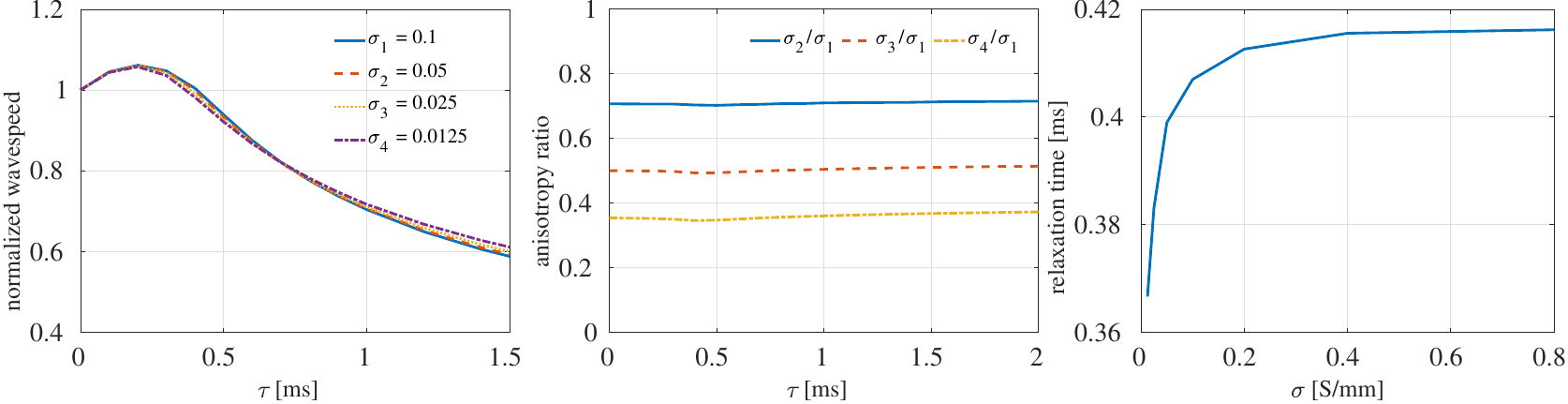}
\caption{Variation of the conduction velocity with respect to the conducitivity coefficients using  the Fenton-Karma model and parameter set 3 in Table \ref{Tab:FK}. Left)  wave speed for different conductivities, normalized with respect to the parabolic conduction velocity. Center) Ratio between the conduction velocities as a function of the relaxation time. The variations are within 5\% 
of the mean. Right) Relaxation time at which the hyperbolic conduction velocity matches  the parabolic conduction velocity. 
Although the relaxation times at which we find the same wave speed for the parabolic and hyperbolic monodomain models are different for different conductivities, the error in the conduction velocity we commit considering the fixed value $\tau = 0.4$ ms  is smaller than 2\%
and therefore typically smaller than the spatial error.
}
\label{fig:test3}
\end{figure*}

\subsection{Conduction velocity anisotropy\label{sub:resultsanisotropy}}

This section analyzes whether inductances affect the anisotropy in conduction velocity.
We define the anisotropy ratio as the ratio between the conduction velocities in the transverse and longitudinal directions, and to obtain a simple estimate for this ratio, we consider plane-wave solutions propagating in these directions.
We use the Fenton-Karma model with parameter set 3 and consider four values of the conductivity: $\sigma_1 = 0.1$ mS/mm, $\sigma_2 = 0.05$ mS/mm, $\sigma_3 = 0.025$ mS/mm, and $\sigma_4 = 0.0125$ mS/mm. Evaluating the wave speeds $v_1$, $v_2$, $v_3$, and $v_4$ corresponding to the  conductivities $\sigma_1$, $\sigma_2$, $\sigma_3$, and $\sigma_4$, respectively, we show how the anisotropy ratio changes with respect to the relaxation time for several conductivity ratios.  
In fact, as shown in 
Table \ref{Tab:FK}, $\sigma_1$ and $\sigma_4$ are the typical longitudinal and transverse conductivities for this model.  
Considering the interval $\Omega=\left[0,5\right]$
cm, we use a mesh size  $h=25$ $\mu$m and a time step size
$\Delta t=0.0025$ ms. 
Fig.~\ref{fig:test3}(left) shows how the velocities at different $\sigma$ compare to each other.
To make the comparison more clear, we normalize the results with respect to the conduction velocities of the parabolic monodomain model, i.e., corresponding to $\tau = 0$ ms.
The results shown in Fig.~\ref{fig:test3}(left) should be understood as the percentage difference of the wave speed with respect to the parabolic case.
It is clear that the curves are similar for the conductivities considered. On the other hand, the relaxation time at which the conduction velocity of the hyperbolic monodomain model matches the conduction velocity of the parabolic monodomain model decreases with smaller conductivities.
Fig.~\ref{fig:test3}(right) shows the relaxation time needed to maintain the same conduction velocity as in the parabolic monodomain model.
In particular, for $\sigma_4$, we find that the relaxation time needed in the hyperbolic monodomain model to yield the same velocities as in the parabolic monodomain models is about $\tau = 0.38$ ms.
The differences in the velocities computed with $\tau =0.4$ ms and $\tau = 0.38$ ms are smaller than 2\%.
This difference is typically smaller than the numerical error in the simulations.
For this reason, in the following tests, our comparison will only use the value $\tau = 0.4$ ms.
Fig.~\ref{fig:test3}(center) shows how the anisotropy ratio is influenced by the relaxation time.  Although the ratio is not constant, the variation over the relaxation times considered never exceeds 5\%.
For the conductivity ratios 8:1 ($\sigma_1$:$\sigma_4$), 4:1 ($\sigma_1$:$\sigma_3$), and 2:1 ($\sigma_1$:$\sigma_2$), we find that the ratios between the conduction velocities in the transverse and in the longitudinal directions are approximately  1/$\sqrt{8}$ ($v_4$/$v_1$), 1/$2$ ($v_3$/$v_1$), and 1/$\sqrt{2}$ ($v_2$/$v_1$).  This behavior is expected because the conduction velocity depends on the square root of the conductivity. We conclude that the effect of the inductances on the anisotropy ratio is negligible.

\subsection{Discretization error: a simple two-dimensional benchmark\label{sub:results_error}}

It is important to be aware of the limitations of the numerical method
used in the simulations. For example, spiral break up can occur easily if the wavefront
is not accurately captured. In fact, numerical error can introduce spatial inhomogeneities
that lead both to spiral wave formation and also to spiral break up.
This numerical artifact disappears under grid refinement, and on sufficiently fine computational meshes, spiral wave break up results only from physical
effects, such as conduction block. 

\begin{figure}[t]
\includegraphics{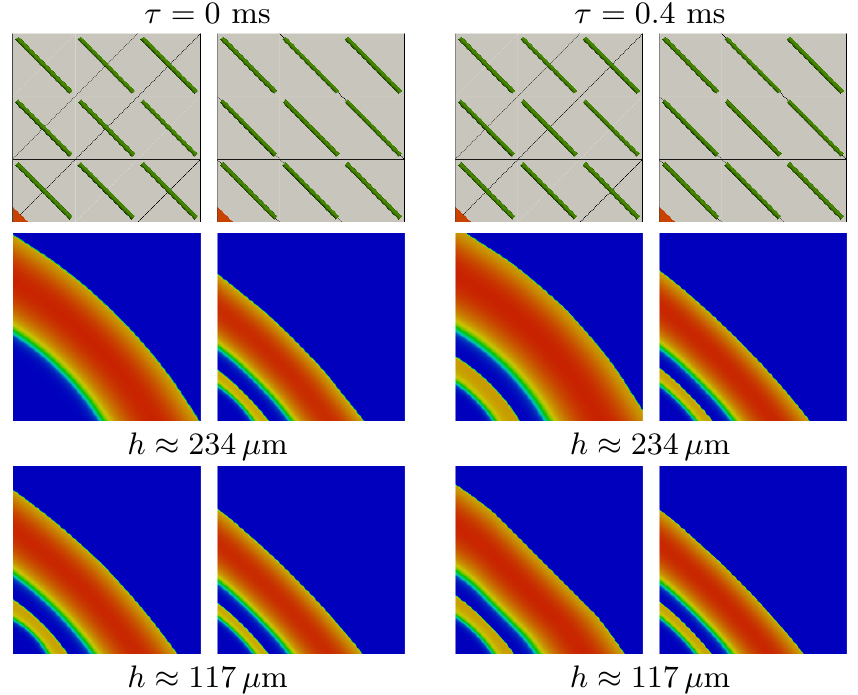}
\protect\caption{Effect of fiber direction and mesh orientation on propagation. Top)
fibers field (green), mesh orientation, and stimulus region (red);
Left) Monodomain model; Right) Hyperbolic monodomain model. Domain
size is 12x12 cm and the solution is evaluated at time $t=450$ ms.}
\label{Fig:meshtest1}
\end{figure}

\begin{figure}[t]
\includegraphics{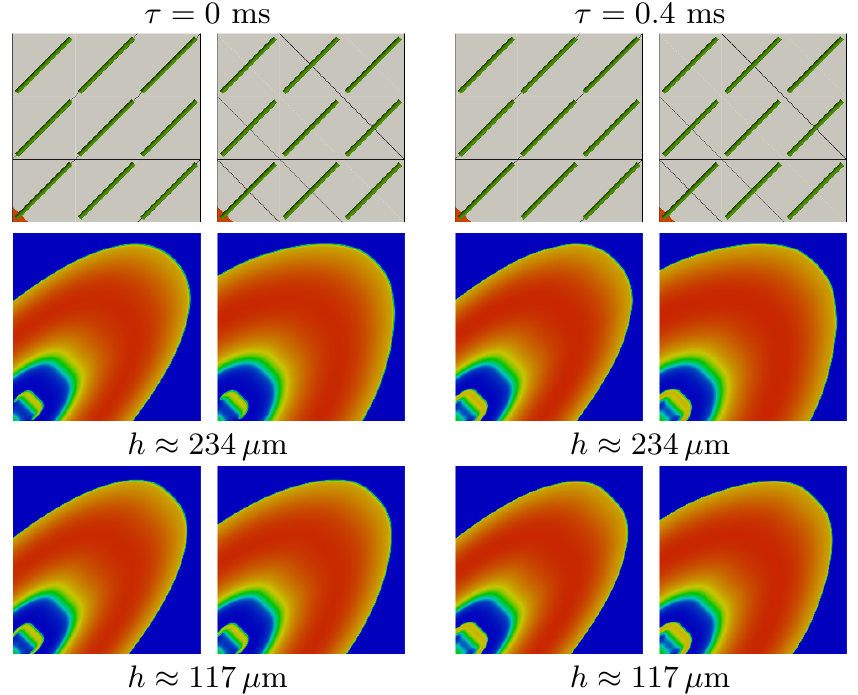}
\protect\caption{Effect of fiber direction and mesh orientation on propagation. Top)
fibers field (green), mesh orientation, and stimulus region (red);
Left) Monodomain model; Right) Hyperbolic monodomain model. Domain
size is 12x12 cm and the solution is evaluated at time $t=360$ ms.}
\label{Fig:meshtest2}
\end{figure}

To examine the role of spatial discretization on the system dynamics, we consider a square slab of tissue, $\Omega=\left[0,12\right]\times\left[0,12\right]$
cm, and we use the Fenton-Karma model with the parameter set 3.
  Reentry is induced by an S1-S2 protocol.  First, an external stimulus of unitary magnitude is applied in
the left bottom corner, $\Omega_{\text{stim}}=\left\{ \boldsymbol{x}\in\Omega:\left\Vert \boldsymbol{x}\right\Vert _{1}\leq1\right\} $
at $t=0$ ms. A second stimulus with the same amplitude and in the
same region is applied after 300 ms. We consider two cases: 1) the fiber
field is aligned with some of the mesh edges; and 2) the fiber field is
not aligned with any edge in the mesh. Figs.~\ref{Fig:meshtest1}(top)
and \ref{Fig:meshtest2}(top) demonstrate that this can be easily achieved by fixing
the fiber field and rotating the mesh. Those figures show the fiber field in green on
top of the computational grid. The red region at the bottom left of the domain is
the stimulus region, $\Omega_{\text{stim}}.$ In the first test
(Fig.~\ref{Fig:meshtest1}), the fibers are set to be orthogonal to the propagating
front. Using 512 elements per side, which is equivalent to a mesh size $h$
of approximately 234 $\mu$m (163 $\mu$m in the direction of wave
propagation), the solutions obtained on the two meshes differ substantially.
Large differences can also be found when we use 1024 elements per side, which is
equivalent to a mesh size $h$ of approximately 117 $\mu$m (82.5 $\mu$m in the
direction of wave propagation). It is clear from Fig.~\ref{Fig:meshtest1}
that whenever the fiber field is not aligned with the mesh, the conduction
velocity is largely overestimated. The second test is similar to the
first one, but the fibers are now rotated by 90$^{\circ}$; see Fig.~\ref{Fig:meshtest2}.
It is clear that the error on the conduction velocity
in the fiber direction is much smaller than the error in the transverse
direction. In fact, the solutions for $h\approx\text{234 }\mu\text{m}$ and $h\approx\text{117 }\mu\text{m}$
show the greatest differences for transverse propagation.

Although similar results have already been reported in prior work
\cite{Krishnamoorthi2013,Niederer2011,Pezzuto2016}, it is nonetheless generally
believed that a mesh size of the order of 200 $\mu$m is usually sufficient
for cardiac computational electrophysiology\cite{Krishnamoorthi2014,fastl2016personalized}. 
This belief is supported
by numerical simulations using a mesh size of 200 $\mu$m that show that the error on the wave speed
in the fiber direction is less than 5\%. On the other hand, the
most interesting dynamics of the propagating front often will occur perpendicular
to the fiber direction. For example, during the normal electrical
activation of the ventricles, the signal spreads from the endocardium
to the epicardium traveling across the ventricular wall, perpendicular
to the fibers. The reduced conductivity in the transversal direction,
usually taken to be about 8 times smaller, requires a finer resolution
of the grid. As noted by \citet{Quarteroni2017}, the required grid
resolution to capture the transverse conduction velocity with an error
smaller than 5\% is about 25 $\mu$m. These results strongly indicate that to capture
correctly the wave speeds, the mesh discretization needs to be about
the size of a single cardiomyocyte, as also discussed by \citet{hand2010adaptive}.

\begin{figure}[htb]
\includegraphics{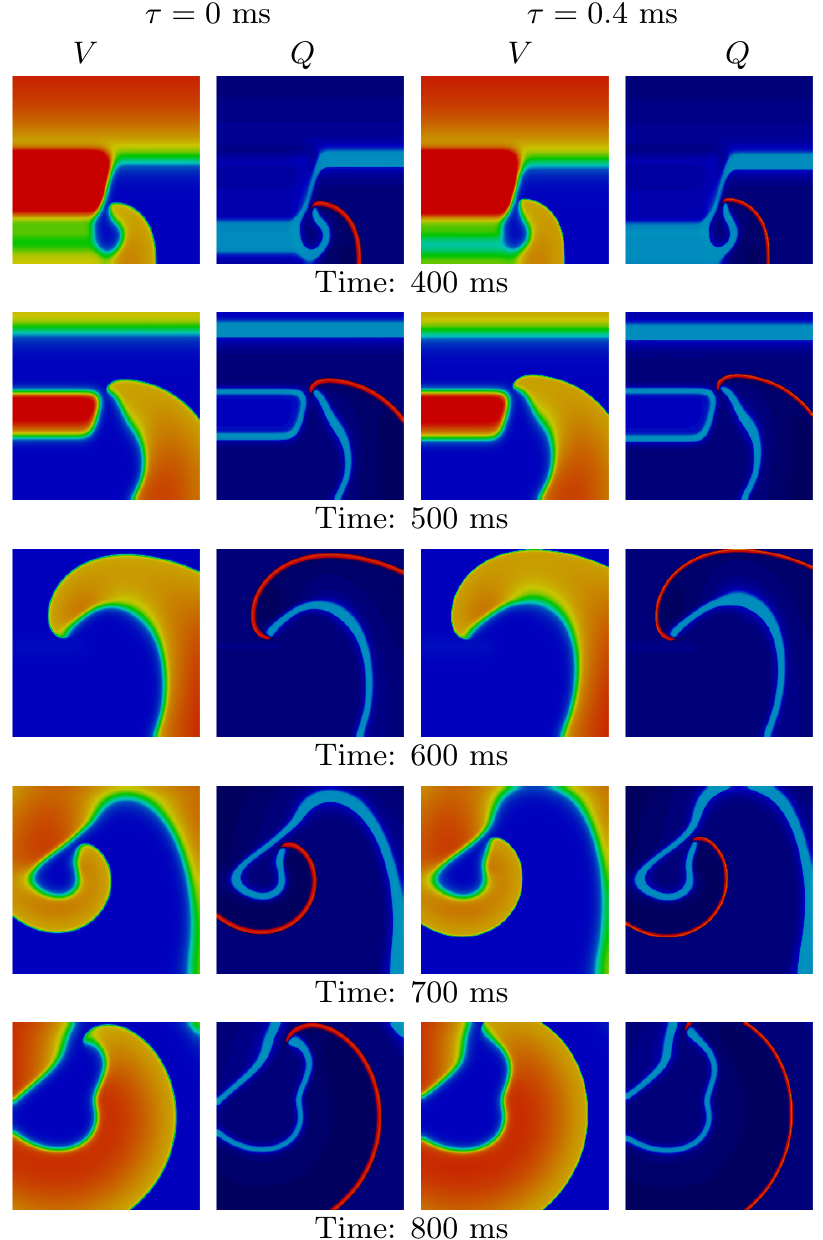}
\protect\caption{Comparison of the evolution of a spiral wave for the isotropic monodomain
model (left) and the isotropic hyperbolic monodomain model (right) using the
Fenton-Karma model with parameter set 3. The wavefronts (red) and
the tails (light blue) represented by the variable $Q$ are shown
next to the transmembrane potential $V$. 
The evolution of the spiral wave  is similar in the two cases. The domain is $\left[0,12\right]\times\left[0,12\right]$
cm, the mesh size is $h\approx100$ $\mu$m, and the time step size $\Delta t = 0.125$ ms. The fibers are aligned
with the $y$-axis.}
\label{Fig:isospiral}
\end{figure}

\begin{figure}[htb]
\includegraphics{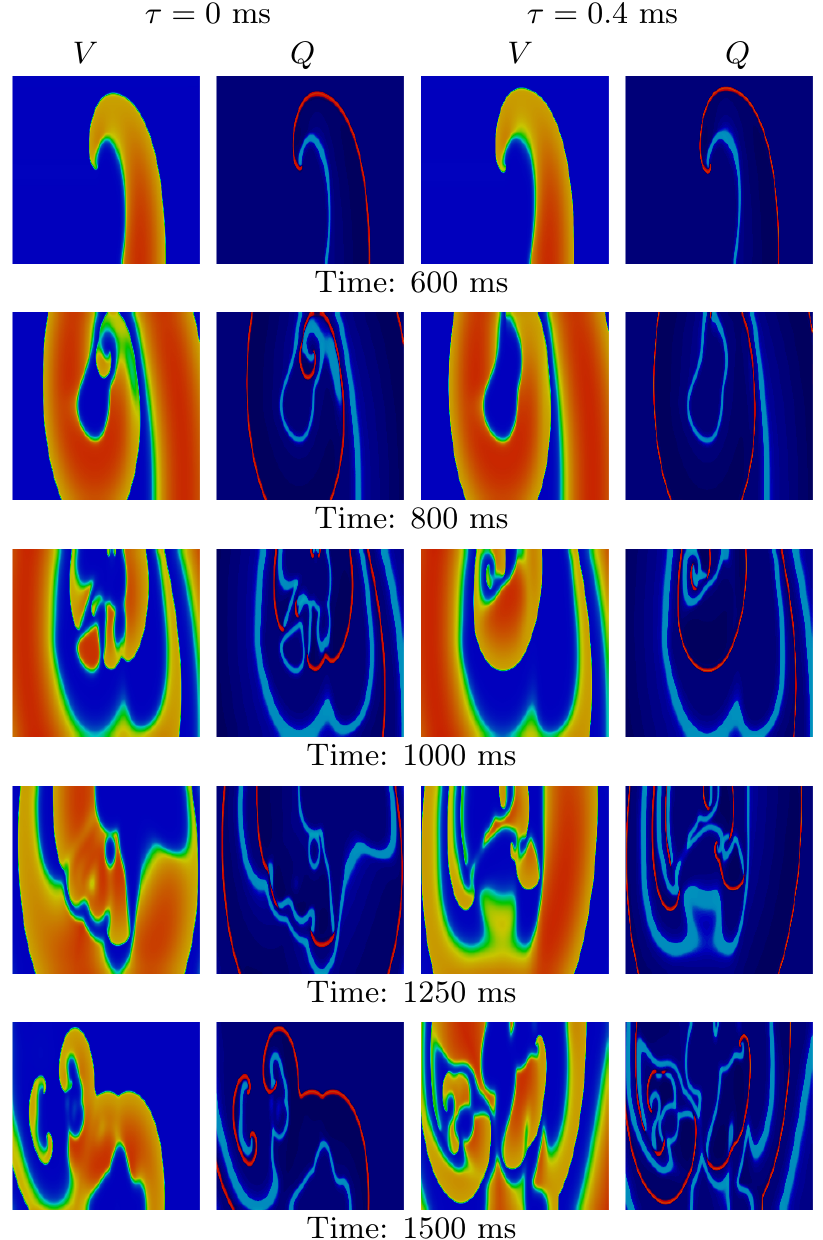}
\protect\caption{Comparison of the evolution of a spiral break up for the monodomain
model (left) and the hyperbolic monodomain model (right) using the
Fenton-Karma model with parameter set 3. The wavefronts (red) and
the tails (light blue) represented by the variable $Q$ are shown
next to the transmembrane potential $V$. Although the initiation
sequence is the same and the spiral wave at $t=600$ ms are similar,
the evolution of the break up is different. The domain is $\left[0,12\right]\times\left[0,12\right]$
cm, the mesh size is $h\approx200$ $\mu$m in the fiber direction, and $h\approx50$ $\mu$m in the transverse direction, and the time step is $\Delta t = 0.125$ ms. The fibers are aligned
with the $y$-axis.}
\label{Fig:spiral}
\end{figure}


%
%
%
%
%
%
%
%

\begin{figure*}[ht]
\includegraphics{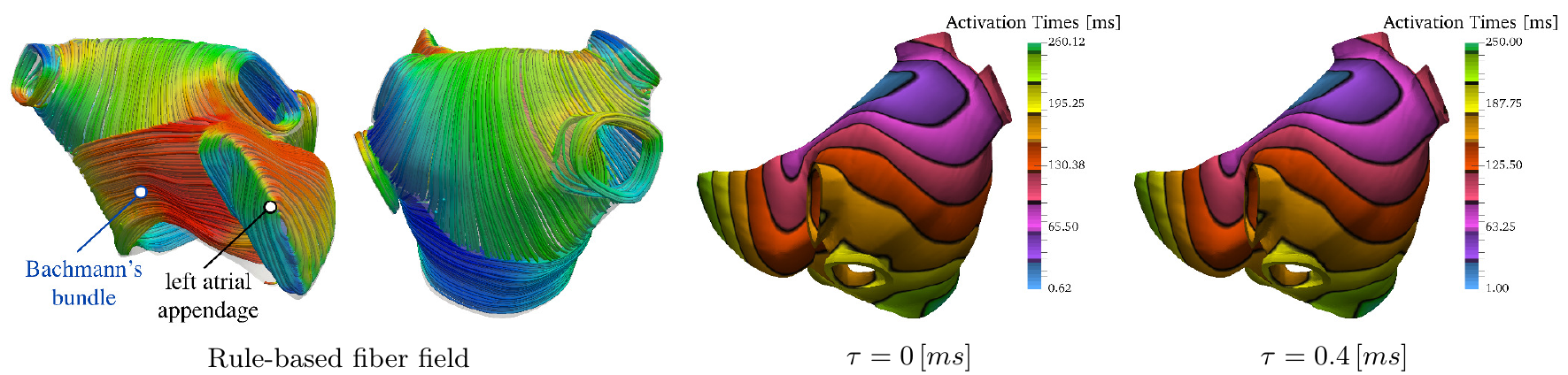}
\protect\caption{Left) Reconstructed fiber field of left atria based on the data in
. Right) Activation times for the monodomain and for the hyperbolic
monodomain model with $\tau=0.4$ ms using the Fenton-Karma model
with parameter set 3.}
\label{Fig:atrialfibersat}
\end{figure*}


\subsection{Effect of the relaxation time on spiral break up\label{sub:results_spirals}}

Our next tests explore the effect of relaxation time on spiral wave break up.
We consider a square slab of tissue, $\Omega=\left[0,12\right]\times\left[0,12\right]$
cm. An initial stimulus
is applied in the region $y<0.5$ cm for 1 ms at $t=0$ ms to generate
a wave propagating in the $y$-axis. A second stimulus is then
applied in the region $\left\{ x<6\text{ cm}\,\wedge\,y<7\text{ cm}\right\} $
for 1 ms at time $t=320$ ms to initiate a spiral wave.
Spiral break up is easily obtained using the parameter set 3 in Table
\ref{Tab:FK}, because of the steep action potential duration (APD) restitution curve. In this case, the back of the wave forms scallops,   and when the turning spiral tries to invade these regions, it encounters refractory tissue and
breaks. 

Fig.~\ref{Fig:isospiral} shows the formation of the spiral wave in the case of isotropy.
With $\tau =0.4$ ms and the Fenton-Karma model with parameter set 3, the conduction velocity of the hyperbolic monodomain model is very close to the conduction velocity of the parabolic monodomain model. In fact, the evolution of the spiral waves in the two cases is very similar.

Fig.~\ref{Fig:spiral} compares spiral break up
obtained using the monodomain model and the hyperbolic monodomain
model in the anisotropic case, with fibers aligned with the $y$-axis.
The relaxation time for the hyperbolic monodomain model is
chosen to be $\tau=0.4$ ms. As explained in Sec.~\ref{sub:resultsanisotropy}, although the conduction velocities in the transverse direction are different in the hyperbolic and parabolic models for $\tau=0.4$ ms, the difference is smaller than 2\%.
This difference is smaller than the spatial error and, for this reason, we consider the two cases to yield essentially the same anisotropy ratio.
It is clearly shown in Fig.~\ref{Fig:spiral} that in the initial phase (up to $t=600$ ms),
the spirals are almost identical. After the first rotation, however, the spiral
wave starts to break, and the relaxation time of the hyperbolic monodomain
model shows quantitative differences in the form of the break up.
Fig.~\ref{Fig:spiral} shows the dynamics of
$Q=\partial_{\text{t}}V$. This variable can be used to define the fronts
(red) and the tails (light blue) of the waves. 

\begin{figure*}[ht]
\includegraphics{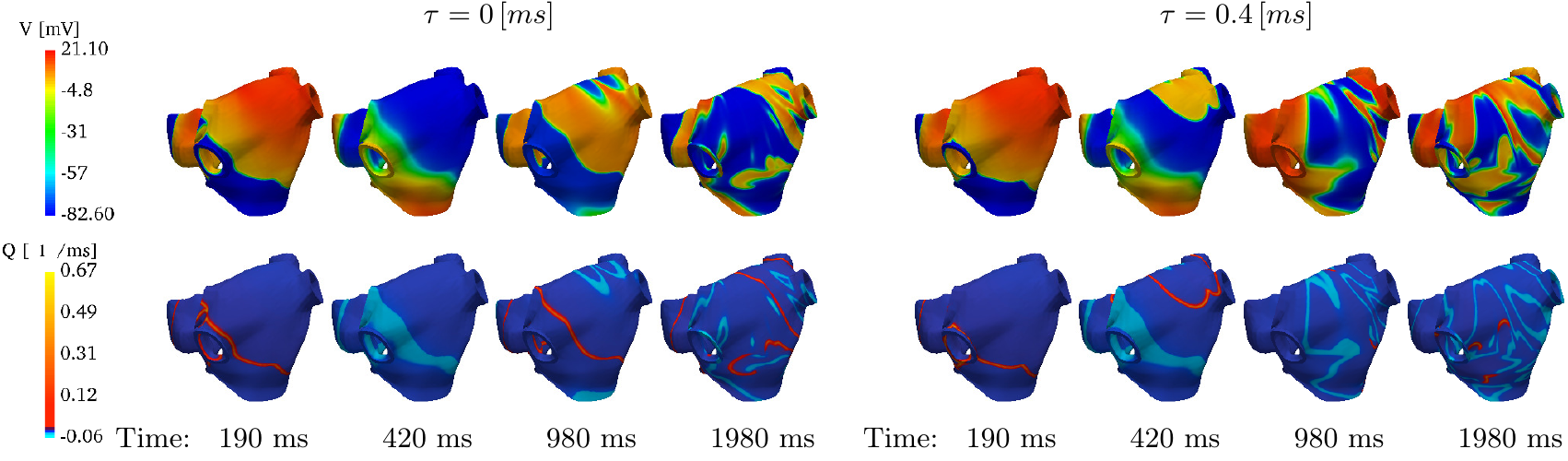}
\protect\caption{Comparison of spiral break up in the left atria using a fast pacing
S1-S2 stimulation protocol using the Fenton-Karma model with parameter
set 3. On the left we show the solution of the monodomain model, relaxation
time $\tau=0$ ms, and on the right we show the solution of the hyperbolic
monodomain model, relaxation time $\tau=0.4$ ms. Using the hyperbolic
monodomain model spiral waves form after second stimulus. The parabolic
monodomain model requires 3 stimuli before spiral waves are formed.
The top row show the solution for the transmembrane potential, and
the bottom row shows the wavefronts (red) and the tails (light blue)
represented by the variable $Q$.}
\label{Fig:atriaspirals}
\end{figure*}


\subsection{Atrial fibrillation\label{sub:resultsatria}}

As a more complex application of the model, we use the hyperbolic monodomain model to simulate
atrial fibrillation. Although a rigorous study of atrial fibrillation
requires the use of a bidomain model, similar to the one proposed in
equations (\ref{eq:HHbidomainH1}) and (\ref{eq:HHbidomainH2}),
the results given by the monodomain model can represent a reasonable
approximation to the dynamics of the full bidomain model. For instance, \citet{Potse2006} provide a comparison of the dynamics of the bidomain and monodomain
models.

The anatomical geometry used in these simulations was based on the human heart model constructed by \citet{segars20104d} using a 4D extended
cardiac-torso (XCAT) phantom. From their data, we extracted and reconstructed
a geometrical representation of the left atrium using SOLIDWORKS. We then
used Trelis to generate a simplex mesh of the left atrium consisting
of approximately 3.5 million elements. 
Our simulations use the Fenton-Karma model with parameter set 3.

\begin{figure*}[ht]
\includegraphics{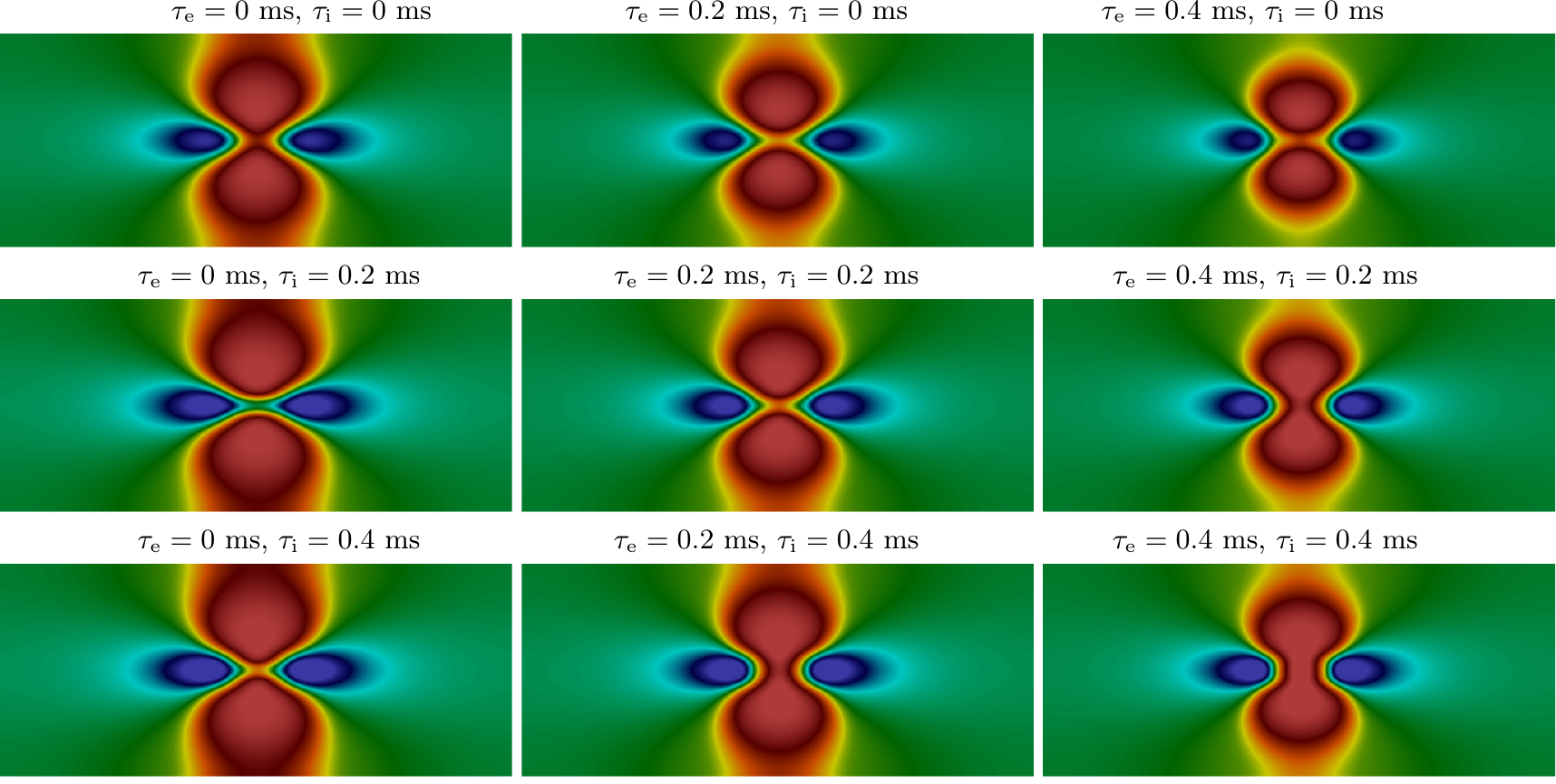}
\protect\caption{Transmembrane potential patterns formed by the
hyperbolic bidomain model with the Fenton-Karma model at 2 ms after a unipolar cathode stimulus is supplied in the
extracellular compartment.
The red regions are areas of depolarization and blue regions are areas of hyperpolarization. The hyperbolic bidomain
model shows a behavior qualitatively similar to the conventional bidomain model shown in the study by \citep{wikswo2009virtual}, but
the details of the are influenced by the extracellular and intracellular relaxation times.
}
\label{Fig:VEP}
\end{figure*}

A major challenge in modeling the atria
is the definition of the fiber field. In work spanning the past 100 years, several studies have tried to characterize
the fiber architecture of the atria \cite{papez1920heart,pashakhanloo2016myofiber,sanchez2013standardized},
but the structure of the muscle in the atria is so complex (as
can be seen from the diagrams by \citet{thomas1959muscular}) that
developing a set of mathematical rules to reproduce a realistic fiber field is challenging. Our strategy
for modeling the fiber field is to use the gradient of a harmonic function that satisfies prescribed boundary conditions. We divided the left atrium
into several regions, and solved three Poisson problems: one for the
left atrial appendage; one for the Bachmann's bundle; and one for the
remaining parts of the left atrium. 
In particular, for the left atrial
appendage, we used a method similar to the one proposed by \citet{Rossi2014}
In the other regions, similar to work by \citet{Patelli2015} and \citet{krueger2011modeling},
we used the normalized gradient of the solution
of the Poisson problem to reconstruct a fiber field similar to the
one shown in the studies by Ho et al.~\cite{ho2002atrial, ho2012left}. The
approximate fiber field of the left atrium reconstructed using this
strategy is shown in Fig.~\ref{Fig:atrialfibersat}, where the
colors represent the magnitude of the $x$-component of the fibers
and are used only to highlight changes in direction. 

To initiate atrial fibrillation, we use an S1-S2 stimulation
protocol that is applied at the junction with the inter-atrial band, as
described by \citet{Colman2014a}.
The S1 stimulus is applied with a cycle length
of 350 ms followed by a short S2 stimulus with a cycle length of 160 ms.
Fig.~\ref{Fig:atrialfibersat}(right) shows that the
initial electrical activation times are similar for the monodomain
and hyperbolic monodomain models. Because this set of parameters gives a velocity
of about 45 cm/s, roughly 35\% slower than the expected conduction
velocity \cite{harrild2000computer}, it takes longer for the atrium
to be fully activated. When the second stimulus is applied, spiral
waves form in the hyperbolic monodomain system but not in the parabolic
monodomain model; see Fig.~\ref{Fig:atriaspirals} at 420 ms. Spiral
waves are generated for $\tau=0$ ms only after the third stimulus
is applied. Fig.~\ref{Fig:atriaspirals} also shows the fronts
(red) and tails (light blue) of the waves, as indicated by
$Q$.


\subsection{Virtual electrodes \label{sub:VEP}}
  The reduction of the
bidomain model to the monodomain model is not valid if the extracellular and intracellular compartments have
different anisotropy ratios. As shown in experimental
studies\cite{knisley1995transmembrane,wikswo1995virtual}, when a current is supplied by an
extracellular electrode, adjacent areas of depolarization and
hyperpolarization are formed in the case of unequal anisotropy ratios. Applying a stimulation current $I_{\text{e,stim}}$ in the extracellular compartment, the
region of hyperpolarization near a cathode ($I_{\text{e,stim}}<0$) is called a virtual anode, while the region of depolarization near an anode ($I_{\text{e,stim}}>0$) is called a virtal cathode. The existence of virtual cathodes and anodes is important to understand the four mechanisms of cardiac stimulation and to study the
mechanisms of defribillation. In this section, we focus only on a unipolar cathodal stimulation to investigate the cathode formation
mechanism in the hyperbolic bidomain model. 
We refer the reader interested in this
virtual electrode phenomenon to the more detailed numerical studies of \citet{wikswo2009virtual} and \citet{franzone2012cardiac}.

Following \citet{sepulveda1989current}, we consider the domain $\Omega =
[-2,2]~\text{cm}\times[-0.8,0.8]~\text{cm}$, and consider this to correspond to the epicardial surface.
We apply a unipolar cathodal current in the region $\Omega_{\textrm{stim}} =
[-0.5,0.5]~\text{mm}\times[-0.1,0.1]~\text{mm}$.
We use the Fenton-Karma model with parameter set 3, and the corresponding nondimensionalized current stimulus amplitude is $-$100.
We choose $\sigma^{\textrm{f}}_\text{e} = 1.5448$, $\sigma^{\textrm{s}}_\text{e}=\sigma^{\textrm{n}}_\text{e}
=1.0438$ and $\sigma^{\textrm{f}}_\text{i} = 2.3172$, $\sigma^{\textrm{s}}_\text{i}=\sigma^{\textrm{n}}_\text{i} = 0.2435$, with the fibers
aligned with the $x$-axis. 

Fig.~\ref{Fig:VEP} shows the pattern formed by the transmembrane potential 2 ms after the extracellular stimulus is
applied. Red regions are areas of depolarization, and blue regions are areas of hyperpolarization. The depolarization
and hyperpolarization regions are affected by the choice of the extracellular and intracellular relaxation times.  For
this test, we consider the relaxation times $\tau_\text{i}$ and $\tau_\text{e}$ to take the values 0, 0.2, and 0.4 ms.
Note that if $\tau_\text{i} = \tau_\text{e}$, the hyperbolic bidomain equations \eqref{eq:HH2bidomainH1} and
\eqref{eq:HH2bidomainH2} reduce to the hyperbolic-elliptic system of equations
\eqref{eq:HEbidomainH} and \eqref{eq:HEbidomainE}.  Although this test shows the ability of the hyperbolic
bidomain model to reproduce the virtual electrode phenomenon, the influence of the relaxation times in the
stimulation remains unclear and necessitates further investigation.


\section{Conclusions\label{sec:Conclusions}}

Local perturbation in the standard bidomain and monodomain models
propagate with an infinite speed. To correct this unrealistic feature, we developed a hyperbolic
version of the bidomain model, and in the case that intracellular and extracellular conductivity tensors have equal anisotropy ratios, we reduced this model to a hyperbolic monodomain model. Our derivation
relies on a Cattaneo-type model for the fluxes, described
by equations (\ref{eq:cattaneofluxe}) and (\ref{eq:cattaneofluxi}).
The Cattaneo-type fluxes are equivalent to introducing self-inductance
effects, as shown by the schematic diagram in Fig.~\ref{Fig:circuit}. As shown in Appendix~\ref{Appendix:Cable-equation}, relaxation times introduced in the
Cattaneo fluxes represent the ratio between the inductances and
the resistances of the extracellular and intracellular compartments.
Although the hyperbolic monodomain model reduces to the classical parabolic model in the case that the relaxation times of both the intracellular and extracellular compartments are zero, the models differ if at least one of these relaxation times is nonzero. 

Although hyperbolic bidomain and monodomain models do not appear to 
have been previously discussed in the literature, the work of \citet{Hurtado2016} does
address the problem of propagation with infinite speed. Their approach
uses a nonlinear model for the fluxes based on porous medium assumptions.
Because of the nonlinear nature of the fluxes, the porous medium approach
to cardiac electrophysiology is difficult to analyze, even for simple
linear reactions. \citet{Hurtado2016} use the simplified ionic model
proposed by \citet{Aliev1996}. It would be interesting to
see the extension to biophysically detailed models. Moreover, it is clear that
the porous medium approach is not incompatible with the theory developed
here. In fact, the two models could be combined using a nonlinear
Cattaneo-type model for the fluxes. 

As discussed by \citet{King2013Determinants}, abnormal conduction
velocities may contribute to arrhythmogenesis. We have studied how
the conduction velocities of the propagating action potential change
with respect to the relaxation time of the Cattaneo-type fluxes. In
the case of the piecewise-linear model of \citet{McKean1970}, it
is possible to find an analytical expression for the conduction velocity of the hyperbolic
monodomain model. (Unfortunately, the wave speed of this bistable model has
been reported incorrectly in some prior studies \cite{Abdusalam2004,MaEndez1998}.)
We have shown that for the hyperbolic monodomain model, signals
cannot travel faster than the characteristic propagation speed of the
medium. Our numerical results match the theoretical predictions for
different values of the excitability parameter. In this simple case,
conduction velocity is a monotone decreasing function of relaxation time.
Additionally, we have shown that, for this linear case, discretization errors are lower in the hyperbolic model.

Using simple nonlinear ionic models \cite{Aliev1996,Fenton1998}, however, we found that the
relationship between the relaxation time and the conduction velocity
is not always monotone. In fact, for small and moderate values of
the relaxation time, we found that waves propagate faster than in the
standard monodomain model. With these ionic
models, a maximum conduction velocity is found at an optimal relaxation time. 
At small relaxation times, the second derivative terms are relatively small.  Consequently, this phenomenon likely results from
the time derivative of the ionic currents. For larger values of the relaxation times, inertial effects dominate, and the conduction velocity monotonically decreases
to zero. This implies that we can find a value of the relaxation time
for which the action potential propagates at the same velocity as
in the standard monodomain model. This fact allowed us to directly
investigate the influence of inductances in the monodomain model.
In the biophysically detailed ionic models for ventricular myocytes
\cite{TenTusscher2006} and atrial myocytes \cite{Grandi2011} considered here,
a similar behavior can be found, although only for very
small values of the relaxation time. For this reason, we compared
the parabolic and hyperbolic monodomain models
in two and three spatial dimensions using the Fenton-Karma model \cite{Fenton1998}.

Changing the type of equation may have important consequences in terms of the stability of the numerical
approximations. 
We did not encounter any substantial numerical difficulties beyond those already present in the standard parabolic models. 
On the contrary, we have shown that the linear system is easier to solve and the solutions are more accurate when the
hyperbolic system is used.
Moreover, we believe that the solution of the parabolic monodomain and bidomain model already have many of the numerical
difficulties usually associated with hyperbolic systems: the propagating action potential in the parabolic models
already have sharp fronts which represent one of the main difficulties in hyperbolic systems. 
Additionally, the hyperbolic monodomain and bidomain models are damped wave equations, where the damping term is 
dominant.
In conclusion, the introduction of hyperbolic terms to the monodomain and bidomain models do not appear to add substantial numerical difficulties. 

We performed several qualitative comparisons between the parabolic and
the hyperbolic monodomain models. In particular, we considered
the case of spiral break up resulting from conduction block caused by a steep
APD restitution curve. Using parameter set 3 of the Fenton-Karma model,
the back of the wave easily forms a series of indentations or scallops.
When a spiral wave invades this after turning, it encounters refractory
tissue that causes conduction block, leading to a break up of the
wave. This behavior, already observed for the standard monodomain
model and explained by \citet{Fenton2002}, is still present in the
hyperbolic monodomain model. We showed spiral break up in a simple
two-dimensional test case and also in an anatomically realistic three-dimensional model of the left atrium. We reconstructed
on the left atria a fiber field qualitatively in accordance with the
data reported by \citet{ho2002atrial,krueger2011modeling}, and \citet{pashakhanloo2016myofiber}
We applied a fast-pacing S1-S2 stimulation protocol to the left atrium
to induce spiral waves and spiral break up. Although the initial activation
sequences were similar, for the hyperbolic monodomain model, spiral
waves appeared already after the second stimulus was applied, whereas
for the parabolic monodomain at least three stimuli were needed.

We also used a simple two-dimensional benchmark to investigate how the alignment
of the fiber direction with respect to the elements of the finite element mesh influences
wave propagation. Under mesh refinement, spatial discretization
errors become smaller and smaller, but the common assumption
that a grid resolution of 200 $\mu$m is sufficient to correctly capture
the conduction velocity \cite{Krishnamoorthi2013,Krishnamoorthi2014} seems to be overly optimistic. In our tests, we analyzed
the propagation of a simple wave using an anisotropic conductivity
tensor. Fixing the fiber field and rotating the mesh, we were able
to show that using a mesh size of about 234 $\mu$m (163 $\mu$m
in the direction of wave propagation), the propagation was greatly
influenced by the mesh orientation relative to the fiber alignment.
Even on a finer mesh with a mesh size of about 117 $\mu$m (82.5
$\mu$m in the direction of wave propagation), visible differences between the propagation patterns obtained with the different orientations remained.
These tests highlight
the fact that if the transverse conductivity coefficients are taken
to be about 8 times smaller than the longitudinal conductivity, as is typically done in practice \cite{roth1997electrical,clayton2011models,hooks2007laminar,benson2010construction},
the mesh size needs to be much smaller than 200 $\mu$m in order to yield resolved dynamics. Our results
are in accordance with the convergence study by \citet{Rossi2014}
for propagation in the transverse direction, in which the authors
showed that the mesh size in the transverse direction needs to be about
25 $\mu$m. Because such mesh resolutions require elements only slightly
larger than the actual cardiac cells, it is natural to ask whether the
use of a continuum model is even appropriate: for a similar
computational cost, it could be possible to construct a discrete model
of the heart. Although the difficulties in representing correctly
the conduction velocities are well documented in the literature \cite{Pathmanathan2011,Krishnamoorthi2013}
and have been thoroughly analyzed by \citet{Pezzuto2016}, the focus
generally seems to have been on the conduction in the fiber direction, as is
clear from the choice of the three-dimensional benchmark test proposed by \citet{Niederer2011}
Other solution strategies for the monodomain and bidomain models using
adaptive mesh refinement \cite{ColliFranzone2006,krause2015towards}
and high-order elements \cite{arthurs2012efficient} have been
proposed, but their application so far appears to be relatively limited.

Finally, we showed that the hyperbolic bidomain model can capture the virtual electrode phenomenon observed in the standard model. 
By stimulating the epicardial surface using a unipolar cathodal current, we have shown how the pattern of the
transmembrane potential is influenced by the intracellular and extracellular relaxation times.
Nonetheless, a more detailed study on how the relaxation times affect the virtual electrodes is needed to fully
understand the hyperbolic bidomain model. 
Comparing the stimulation patterns to experimental data could even serve to calibrate the magnitude of the relaxation times in the hyperbolic bidomain model.

Despite significant progress in models of cardiac electrophysiology,
the complex effects of nonlinearity and heterogeneity in the
heart are far from being fully understood. We have shown that the nonlinearities
of the underlying physics can give unexpected results, in contradiction
to the linear case. Inductances in tissue propagation could
play an important role in cardiac electrical dynamics, especially close to the wavefronts, where the
fast currents responsible for the initiation of the action potential
can give a small but non-negligible contribution. 
The main phenomenon we observed in the hyperbolic model is the enhancement of the conduction velocity at
small relaxation times because of the presence of the time derivative of the ionic currents.
This phenomenon would allow electric signals to propagate at the same speed with lower conductance as compared to the standard models. Future experimental
studies are needed to confirm the importance of these effects. 

\section*{Acknowledgements}

The authors gratefully acknowledge research support from NIH award HL117063 and from NSF award ACI 1450327. The human
atrial model used in this work was kindly provided by Prof.~W.~P.~Segars.
We also thank Prof.~C.~S.~Henriquez and Prof.~J.~P.~Hummel for extended discussions on atrial modeling.
We also gratefully acknowledge support by the libMesh developers in aiding in the development the code used in the
simulations reported herein.

\appendix

\section{Cable equation with inductances\label{Appendix:Cable-equation}}

We follow the derivation of the cable equation by \citet{Keener2009}.
As shown in Fig.~\ref{Fig:circuit}, we assume there exist inductance
effects in both the extra- and intracellular axial currents. We have
\begin{eqnarray}
V_{\text{e}}\left(x\right)-V_{\text{e}}\left(x+\Delta x\right) & = & -L_{\text{e}}\partial_{t}I_{\text{e}}\Delta x-R_{\text{e}}I_{\text{e}}\Delta x,\label{eq:Ohmse}\\
V_{\text{i}}\left(x\right)-V_{\text{i}}\left(x+\Delta x\right) & = & -L_{\text{i}}\partial_{t}I_{\text{i}}\Delta x-R_{\text{i}}I_{\text{i}}\Delta x,\label{eq:Ohmsi}
\end{eqnarray}
where $I_{\text{e}}$ and $I_{\text{i}}$ are the extracellular and intracellular
axial currents, respectively. Dividing by $\Delta x$ and taking
the limit as $\Delta x\to0$, we find
\begin{eqnarray}
\partial_{x}V_{\text{e}} & = & -L_{\text{e}}\partial_{t}I_{\text{e}}-R_{\text{e}}I_{\text{e}},\label{eq:KVLe}\\
\partial_{x}V_{\text{i}} & = & -L_{\text{i}}\partial_{t}I_{\text{i}}-R_{\text{i}}I_{\text{i}}.\label{eq:KVLi}
\end{eqnarray}
The minus sign on the right hand sides is a convention that ensures that positive
charges flow from left to right. Applying Kirchhoff's current law,
we have that 
\begin{eqnarray}
I_{\text{e}}\left(x\right)-I_{\text{e}}\left(x+\Delta x\right)+I_{\text{t}}\Delta x & = & 0,\label{eq:CLdiscretee}\\
I_{\text{i}}\left(x\right)-I_{\text{e}}\left(x+\Delta x\right)-I_{\text{t}}\Delta x & = & 0,\label{eq:CLdiscretei}
\end{eqnarray}
which, in the limit $\Delta x\to0$, becomes 
\begin{equation}
I_{\text{t}}=\partial_{x}I_{\text{e}}=-\partial_{x}I_{\text{i}}.\label{eq:KCL}
\end{equation}
Differentiating  equations (\ref{eq:KVLe}) and
(\ref{eq:KVLi}) with respect to $x$, and equation (\ref{eq:KCL}) with respect to $t$,
we find that for the extracellular compartment,
\begin{eqnarray}
\partial_{xx}V_{\text{e}} & = & -L_{\text{e}}\partial_{xt}I_{\text{e}}-R_{\text{e}}\partial_{x}I_{\text{e}},\label{eq:Veqe}\\
\partial_{xt}I_{\text{e}} & = & \partial_{t}I_{\text{t}},\label{eq:Ieqe}
\end{eqnarray}
and for the intracellular compartment,
\begin{eqnarray}
\partial_{xx}V_{\text{i}} & = & -L_{\text{i}}\partial_{xt}I_{\text{i}}-R_{\text{i}}\partial_{x}I_{\text{i}},\label{eq:Veqi}\\
\partial_{xt}I_{\text{e}} & = & -\partial_{t}I_{\text{t}}.\label{eq:Ieqi}
\end{eqnarray}
Substituting the mixed derivative of the extracellular and intracellular
currents, we find
\begin{equation}
\partial_{xx}V_{\text{e}}=-L_{\text{e}}\partial_{t}I_{\text{t}}-R_{\text{e}}I_{\text{t}},\quad\partial_{xx}V_{\text{i}}=L_{\text{i}}\partial_{t}I_{\text{t}}+R_{\text{i}}I_{\text{t}}.\label{eq:Veqsie}
\end{equation}
Defining $V=V_{\text{i}}-V_{\text{e}}$, and taking the difference between the intracellular
and the extracellular equations, we can write a single equation for
$V$,
\begin{equation}
\partial_{xx}V=\left(L_{\text{i}}+L_{\text{e}}\right)\partial_{t}I_{\text{t}}+\left(R_{\text{i}}+R_{\text{e}}\right)I_{\text{t}}.\label{eq:Veq}
\end{equation}
Recalling that 
\begin{equation}
I_{\text{t}}=\chi\left(C_{\text{m}}\partial_{t}V+I_{\text{ion}}\right),\label{eq:transmembranecurrent}
\end{equation}
we finally find the hyperbolic form of the cable equation,
\[
\tau C_{\text{m}}\dfrac{\partial^{2}V}{\partial t^{2}}+C_{\text{m}}\dfrac{\partial V}{\partial t}-\dfrac{\partial}{\partial x}\left(D\dfrac{\partial V}{\partial x}\right)=-I_{\text{ion}}-\tau\dfrac{\partial I_{\text{ion}}}{\partial t},
\]
where we define the conductivity as $D=1/\chi\left(R_{\text{i}}+R_{\text{e}}\right)$
and the relaxation time by
\begin{equation}
\tau=\dfrac{L_{\text{i}}+L_{\text{e}}}{R_{\text{i}}+R_{\text{e}}}.\label{eq:relaxtime}
\end{equation}

\section{Exact solution for the piecewise-linear bistable model\label{Appendix:exactsolution}}

Consider the simplified piecewise-linear bistable model, 
\begin{eqnarray}
I_{\text{ion}}\left(V\right) & = & kV\nonumber \\
 & - & k\left[V_{2}H(V-V_{1})+V_{0}\left(1-H(V-V_{1})\right)\right].\qquad\label{eq:MCKeanDimensional}
\end{eqnarray}
The nondimensional form of equation (\ref{eq:MCKeanDimensional})
is the simplified model proposed by \citet{McKean1970}, 
\begin{equation}
\hat{I}_{\text{ion}}\left(\hat{V}\right)=\hat{V}-H\left(\hat{V}-\alpha\right).\label{eq:MCKeanAppendix}
\end{equation}
Using (\ref{eq:MCKeanDimensional}) in the hyperbolic monodomain
model (\ref{eq:Hmonodomain}) and considering only one spatial dimension,
we have 
\begin{equation}
\begin{cases}
\mu\dfrac{\partial^{2}\hat{V}}{\partial\hat{t}^{2}}+\left(1+\mu\right)\dfrac{\partial\hat{V}}{\partial\hat{t}}-\dfrac{\partial^{2}\hat{V}}{\partial\hat{x}^{2}}+\hat{V}=0, & \hat{V}<\alpha,\\
\mu\dfrac{\partial^{2}\hat{V}}{\partial\hat{t}^{2}}+\left(1+\mu\right)\dfrac{\partial\hat{V}}{\partial\hat{t}}-\dfrac{\partial^{2}\hat{V}}{\partial\hat{x}^{2}}+\hat{V}=1, & \hat{V}>\alpha,
\end{cases}\label{eq:nondimensionalmckean}
\end{equation}
 where we have introduced the nondimensional variables
\[
\hat{t}=\dfrac{1}{T}\hat{t},\quad\hat{\boldsymbol{x}}=\dfrac{1}{L}\boldsymbol{x},\quad\hat{V}=\dfrac{V-V_{0}}{V_{2}-V_{0}}.
\]
In particular, we define
\[
T=\dfrac{C_{\text{m}}}{k},\quad L=\sqrt{\dfrac{\sigma}{k\chi}},\quad\mu=\dfrac{\tau k}{C_{\text{m}}},
\]
where $\mu$ is a nondimensional number that characterizes the effect of the inductances in the system. Specifically, $\mu$
is  the ratio between the relaxation time of the system and the characteristic time of the reactions: the larger $\mu$, the more important the effects of the inductances become.
Introducing the change of
coordinates $z=x-ct,$ such that $U\left(z\right)=U\left(x-ct\right)=\hat{V}\left(x,t\right)$, 
the system (\ref{eq:nondimensionalmckean}) is transformed into
\begin{equation}
\begin{cases}
\left(c^{2}\mu-1\right)U_{zz}-c\left(1+\mu\right)U_{z}+U=0, & z>0,\\
\left(c^{2}\mu-1\right)U_{zz}-c\left(1+\mu\right)U_{z}+U=1, & z<0.
\end{cases}\label{eq:xisystemmckean}
\end{equation}
Consider first the case $z>0$. Using $\gamma=c^{2}\mu-1$ and $\beta=-c\left(1+\mu\right)$,
this equation reads 
\[
\gamma V''+\beta V'+V=0,
\]
for which the roots of the characteristic polynomial are
\begin{equation}
\lambda_{\pm}=\dfrac{-\beta\pm\sqrt{\beta^{2}-4\gamma}}{2\gamma},\label{eq:charactpolynroots}
\end{equation}
so that the solution is $U\left(z\right)=A_{+}e^{\lambda_{+}z}+A_{-}e^{\lambda_{-}z}$.
It is easy to verify that the case $z<0$ has solution $U\left(z\right)=B_{+}e^{\lambda_{+}z}+B_{-}e^{\lambda_{-}z}+1$.
The global solution, therefore, is
\begin{equation}
U\left(z\right)=\begin{cases}
A_{+}e^{\lambda_{+}z}+A_{-}e^{\lambda_{-}z}, & z>0,\\
B_{+}e^{\lambda_{+}z}+B_{-}e^{\lambda_{-}z}+1, & z<0.
\end{cases}\label{eq:solutiongeneralmckean}
\end{equation}
For (\ref{eq:solutiongeneralmckean}) to represent a propagating
front, it is necessary that the solution is real and bounded for
any value of $z$. This implies that either $\lambda_{-}<0<\lambda_{+}$
or $\lambda_{+}<0<\lambda_{-}$. Requiring $\beta^{2}-4\gamma>0$,
we find the condition $\gamma<\beta^{2}/4$, which is always satisfied
for the model (\ref{eq:xisystemmckean}) because $\beta^{2}-4\gamma=c^{2}\left(1-\mu\right)^{2}+4>0.$
If $c>0$ then $\beta < 0$. Then $\lambda_{-}>0$ and $\lambda_{+}<0$ 
if 
\begin{equation}
c<\sqrt{\dfrac{1}{\mu}}.\label{eq:boundvelocity}
\end{equation}
Imposing $U\left(\infty\right)=0$, we find $A_{-}=0$, and
imposing $U\left(-\infty\right)=1$, we find $B_{+}=0$, so that
\begin{equation}
U\left(z\right)=\begin{cases}
A_{+}e^{\lambda_{+}z}, & z>0,\\
B_{-}e^{\lambda_{-}z}+1, & z\leq0.
\end{cases}\label{eq:infinityconditions}
\end{equation}
To ensure the continuity of the solution at $z=0$, we fix $U\left(0\right)=\alpha$,
so that $A_{+}=\alpha=1+B_{-}$ and
\begin{equation}
U\left(z\right)=\begin{cases}
\alpha e^{\lambda_{+}z}, & z>0,\\
\left(\alpha-1\right)e^{\lambda_{-}z}+1, & z\leq0.
\end{cases}\label{eq:continuituconditions}
\end{equation}
The first derivative of $U$ is
\begin{equation}
U'\left(z\right)=\begin{cases}
\alpha\lambda_{+}e^{\lambda_{+}z}, & z>0,\\
\left(\alpha-1\right)\lambda_{-}e^{\lambda_{-}z}, & z\leq0.
\end{cases}\label{eq:Derivative}
\end{equation}
Imposing the continuity of $U'$ at $z=0$, we find $\alpha\lambda_{+}=\left(\alpha-1\right)\lambda_{-}$
and therefore 
\begin{equation}
\dfrac{\gamma}{\beta^{2}}=\dfrac{\left(\alpha^{2}-\alpha\right)}{\left(2\alpha-1\right)^{2}}.\label{eq:continuityderivative}
\end{equation}
 Using the definitions of $\gamma$ and $\beta$ in (\ref{eq:continuityderivative}),
we finally obtain an expression for the speed $c$ in terms of $\mu$
and $\alpha$,
\begin{equation}
c=\dfrac{\left(1-2\alpha\right)}{\sqrt{\mu+\left(\alpha-\alpha^{2}\right)\left(\mu-1\right)^{2}}}\leq\sqrt{\dfrac{1}{\mu}}.\label{eq:speedmckean}
\end{equation}
Notice that the characteristic propagation speed is bounded from above by $\sqrt{1/\mu}$.

\section{Finite Element Discretization\label{sec:FEM}}

The common practice to solve the coupled nonlinear system of cardiac electrophysiology is to split the monodomain (or
bidomain) model from the ionic model.
Motivated by the work of \citet{munteanu2009scalable}, we will also decouple the two systems. 
For brevity, we show the discretization only for the monodomain model with a first-order time
integrator, because first-order time integrators are still widely used in the community.
The same approach can be used to discretize the hyperbolic bidomain model. 
Given the subinterval $[t^n, t^{n+1}]$, with $t^0=0$, we define two separate subproblems, one for the ionic model and
one for the monodomain equation connected by the initial conditions.
In the first step, we solve the ionic model \eqref{eq:ionicmodel} for $\boldsymbol{w}^{n+1}$, and we compute the ionic
currents using the updated values of the state variables. 
The ionic currents computed in this way are then used to solve the monodomain system  (\ref{eq:Vequation})--(\ref{eq:Qequation}).

As shown in \citet{Spiteri2010}, the most efficient numerical algorithm to solve the stiff ODEs of the ionic model
strongly depends on the model considered.
We use the forward Euler method for the simplified ionic models of Aliev-Panfilov and Fenton-Karma. 
For the ventricular ionic model of \citet{TenTusscher2006}, we use the Rush-Larsen method \cite{rush1978practical}.
The  atrial ionic model of \citet{Grandi2011} has a stricter restriction on the time step size. 
Consequently, we use the Rush-Larsen method in combination to the backward Euler method for linear equations. 
The time step size is defined as $\Delta t = t^{n+1}-t^n$. 
To obtain a second-order time discretization scheme, the ionic model can be solved using the explicit trapezoidal method (i.e., Heun's method), which is a second-order accurate strong stability preserving Runge-Kutta method.

Once the solution of the ionic model has been found, we solve the
monodomain model using a low-order finite element discretization.
Because (\ref{eq:Vequation}) is a simple ordinary differential equation,
the computational costs for solving (\ref{eq:Hmonodomain}) or the
system (\ref{eq:Vequation})--(\ref{eq:Qequation}) is comparable.
Denoting by $\left(v,w\right)_{\Omega}=\int_{\Omega}vw$ the $L^{2}\left(\Omega\right)$
inner product and using the boundary conditions (\ref{eq:BC}), the
Galerkin approximation of equations (\ref{eq:Vequation})--(\ref{eq:Qequation})
is to find $V^{h},\,Q^{h}$ $\in\mathcal{S}^{h}$ such that
\begin{eqnarray}
\left(\dfrac{\partial V^{h}}{\partial t},\psi^{h}\right)_{\Omega} & = & \left(Q^{h},\psi^{h}\right)_{\Omega},\label{eq:Galerkinveq}\\
\left(\tau C_\text{m}\dfrac{\partial Q^{h}}{\partial t}+C_\text{m}Q^{h},\phi^{h}\right)_{\Omega}\nonumber \\
+\left(\boldsymbol{D}\nabla V,\nabla\phi^{h}\right)_{\Omega} & = & \left(I_\text{ion}^{h},\phi^{h}\right)_{\Omega}\nonumber \\
 &  & +\left(\tau\dfrac{\partial I_\text{ion}^{h}}{\partial t},\phi^{h}\right)_{\Omega},\qquad\label{eq:Galerkinqeq}
\end{eqnarray}
for all $\psi^{h},$ $\phi^{h}$ $\in\mathcal{S}^{h}$, where
\begin{equation}
\mathcal{S}^{h}=\left\{ v^{h}\in C^{0}\left(\bar{\Omega}\right):\left.v^{h}\right|_{K}\in\mathcal{P}^{1}\left(K\right),\;\forall K\in\mathcal{T}^{h}\right\} .
\end{equation}
In practice, we look for a continuous solution that is linear in every
simplex element $K$ in the triangulation $\mathcal{T}^{h}$ of $\Omega$.
Introducing the basis of nodal shape functions $\left\{ N_{A}\right\} _{A=1}^{M},$
with $M=\text{dim}(\mathcal{S}^{h})$, the solution fields are discretized
as
\[
V^{h}=\sum_{A=1}^{M}N_{A}\boldsymbol{\mathsf{V}}_{A},\quad Q^{h}=\sum_{A=1}^{M}N_{A}\boldsymbol{\mathsf{Q}}_{A}.
\]
Equations (\ref{eq:Galerkinveq}) and (\ref{eq:Galerkinqeq})
yield the matrix system
\begin{eqnarray}
\dot{\boldsymbol{\mathsf{V}}} & = & \boldsymbol{\mathsf{Q}},\label{eq:matrixformveq}\\
C_\text{m}\boldsymbol{\mathsf{M}}\left(\tau\dot{\boldsymbol{\mathsf{Q}}}+\boldsymbol{\mathsf{Q}}\right)+\boldsymbol{\mathsf{K}}\boldsymbol{\mathsf{V}} & = & \boldsymbol{\mathsf{F}}+\tau\boldsymbol{\mathsf{L}}.\label{eq:matrixformqeq}
\end{eqnarray}
We use the half-lumping scheme proposed by \citet{Pathmanathan2011},
so that
\[
\boldsymbol{\mathsf{F}}=\boldsymbol{\mathsf{M}}\boldsymbol{\mathsf{I}},\quad\boldsymbol{\mathsf{L}}=\boldsymbol{\mathsf{M}}\boldsymbol{\mathsf{J}}
\]
and
\begin{eqnarray}
\dot{\boldsymbol{\mathsf{V}}} & = & \boldsymbol{\mathsf{Q}},\label{eq:matrixformveq-1}\\
C_\text{m}\boldsymbol{\mathsf{M}}_\text{L}\left(\tau\dot{\boldsymbol{\mathsf{Q}}}+\boldsymbol{\mathsf{Q}}\right)+\boldsymbol{\mathsf{K}}\boldsymbol{\mathsf{V}} & = & \boldsymbol{\mathsf{M}}\left(\boldsymbol{\mathsf{I}}+\tau\boldsymbol{\mathsf{J}}\right),\label{eq:matrixformqeq-1}
\end{eqnarray}
where
\[
\left[\mathsf{M}_{AB}\right]=\left(N_{A},N_{B}\right)_{\Omega},\qquad\left[\mathsf{K}_{AB}\right]=\left(\nabla N_{A},\boldsymbol{D}\nabla N_{B}\right)_{\Omega},
\]
$\boldsymbol{\mathsf{I}}$ and $\boldsymbol{\mathsf{J}}$ are the
nodal evaluation of $I_\text{ion}$ and $\partial I_\text{ion}/\partial t$,
and $\boldsymbol{\mathsf{M}}_\text{L}$ is the lumped mass matrix obtained
by row summation of the mass matrix. Using a simple first-order implicit-explicit
time integrator, the fully discrete system of equations reads
\begin{eqnarray}
\boldsymbol{\mathsf{V}}^{n+1} & = & \boldsymbol{\mathsf{V}}^{\text{n}}+\Delta t\boldsymbol{\mathsf{Q}}^{n+1},\quad\;\,\label{eq:discretematrixformveq}\\
C_\text{m}\left(\tau+\Delta t\right)\boldsymbol{\mathsf{M}}_\text{L}\boldsymbol{\mathsf{Q}}^{n+1}+\Delta t\boldsymbol{\mathsf{K}}\boldsymbol{\mathsf{V}}^{n+1} & = & \Delta t\boldsymbol{\mathsf{M}}\left(\boldsymbol{\mathsf{I}}^{*}+\tau\boldsymbol{\mathsf{J}}^{*}\right)\quad\nonumber \\
 & + & \tau C_\text{m}\boldsymbol{\mathsf{M}}_\text{L}\boldsymbol{\mathsf{Q}}^{\text{n}},\label{eq:discretematrixformqeq}
\end{eqnarray}
where $\boldsymbol{\mathsf{I}}^{*}$ and $\boldsymbol{\mathsf{J}}^{*}$
indicate the values of the ionic currents, and of their derivatives,
evaluated using the updated values for $\boldsymbol{w}^{n+1}$ obtained
by solving the ionic model. Introducing (\ref{eq:discretematrixformveq})
in (\ref{eq:discretematrixformqeq}), we arrive at a single system
for $\boldsymbol{\mathsf{Q}}^{n+1}$,
\begin{eqnarray}
\left[C_\text{m}\left(\tau+\Delta t\right)\boldsymbol{\mathsf{M}}_\text{L}+\Delta t^{2}\boldsymbol{\mathsf{K}}\right]\boldsymbol{\mathsf{Q}}^{n+1}=\;\;\qquad\qquad\qquad\qquad\nonumber \\
\tau C_\text{m}\boldsymbol{\mathsf{M}}_\text{L}\boldsymbol{\mathsf{Q}}^{\text{n}}-\Delta t\boldsymbol{\mathsf{K}}\boldsymbol{\mathsf{V}}^{\text{n}}+\Delta t\boldsymbol{\mathsf{M}}\left(\boldsymbol{\mathsf{I}}^{*}+\tau\boldsymbol{\mathsf{J}}^{*}\right).\qquad\;\;\label{eq:Qsystem}
\end{eqnarray}
The time derivative of the ionic currents, $\boldsymbol{\mathsf{J}}^{*}$,
is approximated nodally using the chain rule as 
\begin{eqnarray}
\dfrac{\partial I_\text{ion}}{\partial t}(t^n) & \approx & \dfrac{\partial I_\text{ion}}{\partial
V}\left(V^{\text{n}},\boldsymbol{w}^{n+1}\right)Q^{\text{n}}\nonumber \\
 &  & +\dfrac{\partial
 I_\text{ion}}{\partial\boldsymbol{w}}\left(V^{\text{n}},\boldsymbol{w}^{n+1}\right)\cdot\boldsymbol{g}(V^{\text{n}},\boldsymbol{w}^{n+1}
 ),\quad\label{eq:nodaldion}
\end{eqnarray}
where the function $\boldsymbol{g}(V, \boldsymbol{w})$ is the right hand side of the ionic model system defined in
equation \eqref{eq:ionicmodel}. After solving (\ref{eq:Qsystem}), we update the voltage equation
(\ref{eq:discretematrixformveq}). 

The scheme used to obtain \eqref{eq:discretematrixformveq} and \eqref{eq:discretematrixformqeq} is the simplest
method -- ARS(1,1,1) -- of a series of implicit-explicit (IMEX) Runge-Kutta (RK) algorithms that are popular for hyperbolic systems \cite{ascher1997implicit,boscarino2009class}.
For a second-order time integrator, we use the H-CN(2,2,2) scheme\cite{boscarino2016high}, which combines Heun's method for the explicit part and the Crank-Nicholson method for
the implicit part.

We have not experienced any difference in time step size restriction between the parabolic and hyperbolic models. In fact, we
have used an IMEX-RK method where only the ionic currents and their time derivatives are treated explicitly in both
parabolic and hyperbolic models. The time step size restriction is dictated only by the ionic model. Nonetheless, the time step size
should be chosen to accurately capture the propagation of the action potential. Denoting with $v$ the conduction
velocity and with $h_m$  the smallest element
size, we chose our time step size such that the CFL condition $\Delta t \leq h_m / v $ holds.

From a purely numerical perspective, in the first-order scheme, the dissipative nature of the backward Euler method can
be used to remove numerical oscillations, while the second-order Crank-Nicholson method will keep such
oscillations.
A dissipative second-order time integrator could be used in place of the Crank-Nicholson method.
In any case, the physical dissipation of the hyperbolic monodomain model reduces such artifacts for first-~and second-order schemes.

One of the main difficulties arising from hyperbolic systems is their tendencies to form discontinuities for which 
numerical approximations typically develop spurious oscillations. Although the standard monodomain and bidomain models
are parabolic systems, if the  wavefront is not accurately resolved then the upstroke of the action potential can act as a discontinuity.
For this reason, the numerical methods used for cardiac electrophysiology should be already suitable for the hyperbolic systems considered here.

\bibliographystyle{abbrv}
\bibliographystyle{agsm}
\bibliography{biblio}

\end{document}